\documentclass[12pt,draftcls,onecolumn]{IEEEtran}

\usepackage{amsfonts,amsmath,amscd,graphicx,hyperref}

\newtheorem{thm}{Theorem}

\newtheorem{lem}{Lemma}
\newtheorem{defn}{Definition}

\newcommand{\RR}{{\mathbb R}}

\newcommand{\abs}[1]{\lvert#1\rvert}

\newcommand{\Rset}{{\mathbb R}}

\newcommand{\XX}{{\mathcal X }}
\newcommand{\YY}{{\mathcal Y }}

\newcommand{\UU}{{\mathcal  U}}
\newcommand{\LL}{{\mathcal  L}}

\newcommand{\dv}[2]{{\frac{\partial #1}{\partial #2}}}

\newcommand{\dotex}{{\frac{d}{dt}}}

\newcommand{\vw}{\textbf{$\omega$}}

\newcommand{\va}{{\mathbf a}}

\newcommand{\vG}{{\mathbf A}_{grav}}
\newcommand{\vg}{{\mathbf a_{grav}}}

\newcommand{\vB}{{\mathbf B}}

\begin{document}

\title{Symmetry-preserving observers}

\author{Silv\`{e}re Bonnabel, Philippe Martin and Pierre Rouchon
\thanks{S. Bonnabel, Ph. Martin and P. Rouchon are with Centre Automatique et Syst\`{e}mes,
\'{E}cole des Mines de Paris, 60 boulevard Saint-Michel, 75272 Paris
CEDEX 06, FRANCE {\tt\small silvere.bonnabel@ensmp.fr}, {\tt\small
philippe.martin@ensmp.fr}, {\tt\small pierre.rouchon@ensmp.fr}}
}
\date{7 February 2007}
\maketitle

\begin{abstract}

 This paper presents three non-linear  observers for
 three examples of engineering interest:  a
 non-holonomic car, a chemical reactor, and an inertial navigation system. For each example,
 the design is based on  physical symmetries. This
 motivates  the theoretical development of invariant observers, i.e, symmetry-preserving
 observers. We consider an observer to consist of a copy of the system
 equation and a correction term,
 and we propose a constructive method (based on the Cartan
 moving-frame method) to find all the symmetry-preserving correction
 terms. The construction relies on  an invariant frame (a classical notion) and on
 an invariant output-error, a less standard notion precisely defined here. For each
 example,  the convergence analysis relies on the use of invariant state-errors, a symmetry-preserving way
  to define the estimation error.

\end{abstract}


\section{Introduction}
Symmetries have been used in control theory for feedback design and
optimal control, see for instance~\cite{fagnani-willems-siam93,
grizzle-marcus-ieee85,respondek-tall-scl02,schaft-siam87,martin-et-al-cocv03,Spong-Bullo-ieee2005}
but much less for observer
design~\cite{aghannan-rouchon-ieee03,aghannan-rouchon-cdc02,maithripala2005,mahony-et-al-dcd05,hamel-mahony-icra06}.
In this paper we use symmetries for observer design and we develop a
theory of  invariant observers. This theory is motivated by three
non-linear examples of engineering interest: a non-holonomic car, an
exothermic chemical reactor and a velocity-aided inertial navigation
system. In each case the symmetries have an obvious physical
interpretation. For the first example we propose  a non-linear
observer which converges for any initial condition except one
(theorem \ref{carobs:thm}). For the second, we design a non-linear
globally convergent observer (theorem \ref{chemobs:thm}). For the
third, the observer is locally convergent around any system
trajectory. Moreover the global behavior is independent of the
system trajectory (theorem \ref{navobs:thm}). This theory may be
applied to many other systems such as those treated
in~\cite{creamer-jas96,mahony-et-al-dcd05,hamel-mahony-icra06} where
the invariance relative to the choice of the reference 3D-frame  is
exploited in observer design and convergence analysis.

The theoretical contribution  of the paper is the following: for the
smooth system with state $x$,  input $u$ and output $y$, invariance
under the action of a  Lie group $G$ is defined and corresponds to a
separate action of $G$ on the state-space, on the input-space and on
the output space. Invariance means that the dynamics $\dotex
x=f(x,u)$ and  the output map $y=h(x)$ remain  unchanged by a change
of state, input, and output coordinates corresponding  to the action
of $G$. We define invariance for an asymptotic Luenberger nonlinear
observer under the action $G$ similarly, where the group acts also
on the estimated space and the estimated output in a similar way.
When the group dimension does not exceed the state dimension  we
propose (theorem~\ref{invobs:char:thm}) a constructive design of the
invariant observer. This construction is based on an invariant frame
and an invariant output-error. Such invariant output-errors
(definition~\ref{inverr:defn}) are introduced here for the first
time and can be computed via Cartan's moving frame method
(theorem~\ref{thm:ioe}). We show how to transform a locally
convergent asymptotic observer around an equilibrium point into an
invariant one with the same first order approximation. To deal with
convergence issues, we introduce invariant state-errors. The three
examples show that these state-errors play key role in the
convergence analysis.

The content of this paper  is as follows:  in section~\ref{th:sec},
we define invariant systems and invariant pre-observers. The general
form of an invariant pre-observer is given in
theorem~\ref{invobs:char:thm}: it relies on invariant output errors
and invariant vector fields. Their explicit construction relies on
the moving frame method \cite{olver-book99}, which is summarized in
subsection \ref{preliminaries:sec}. Around an equilibrium, we show
it is always possible to build an invariant observer whose linear
tangent approximation is any linear asymptotic observer of the
Luenberger type. To study the convergence, we define invariant state
error. It is a way of defining the error equation so that it
preserves the symmetries, while the usual $\hat x-x$ does not in
general. It obeys a differential system where only the invariant
part of the system trajectory appears (theorem~\ref{inverrdyn:lem}).
This property reduces the dimension of the convergence problem and
appears to play a crucial role in the examples. In
section~\ref{examples}, we study in detail three physical examples.

A summary of the state of the art results on symmetries of dynamic
systems can be found in the monograph \cite{ortega-et-al-book2000}.
The notion of invariant observer and invariant output error can be
found in \cite{aghannan-rouchon-cdc02,aghannan-thesis}. Other
preliminary results presented in this paper can be found
in~\cite{bonnabel-rouchon-LN05,bonnabel-et-al:ACC06,bonnabel-et-al:CIFA06}.

 \section{Invariant systems, observers and errors}\label{th:sec}

\subsection{Invariant systems and equivariant outputs}

\begin{defn}\label{group:action:def}
Let $G$ be a Lie Group with identity~$e$ and $\Sigma$ an open set
(or more generally a manifold). A \emph{transformation group}
$(\phi_g)_{g\in G}$ on~$\Sigma$ is a smooth map
\[ (g,\xi)\in G\times\Sigma\mapsto\phi_g(\xi)\in\Sigma \]
such that:
\begin{itemize}
 \item $\phi_e(\xi)=\xi$ for all~$\xi$
 \item $\phi_{g_2}\bigl(\phi_{g_1}(\xi)\bigr)=\phi_{g_2g_1}(\xi)$
 for all $g_1,g_2,\xi$.
\end{itemize}
\end{defn}

Notice $\phi_g$ is by construction a diffeomorphism on~$\Sigma$ for
all~$g$. The transformation group is \emph{local} if $\phi_g(\xi)$
is defined only when $g$ lies sufficiently near~$e$. In this case
the transformation law
$\phi_{g_2}\bigl(\phi_{g_1}(\xi)\bigr)=\phi_{g_2g_1}(\xi)$ is
imposed only when it makes sense. All the results of the paper being
local, since based on constant rank assumptions, we consider in this
section  only local transformation groups acting on open sets. When
we say ``for all~$g$" we thus mean ``for all~$g$ sufficiently near
the identity~$e$ of~$G$"; in the same way ``for all~$\xi$" usually
means ``for all generic~$\xi$ in~$\Sigma$". We systematically use
these stylistic shortcuts in order to improve readability.

Consider now the smooth system
\begin{align}
 \label{eq:sys}\dotex x&=f(x,u)\\
 \label{eq:output}y &=h(x,u)
\end{align}
where $x$ belongs to an open subset $\XX\subset\Rset^n$, $u$ to an
open subset $\UU\subset\Rset^m$ and $y$ to an open subset
$\YY\subset\Rset^p$, $p \leq n$.

We assume the signals $u(t),y(t)$ known ($y$ is measured, and $u$ is
measured or known - control input, measured perturbation, constant
parameter).

Consider also the local group of transformations on $\XX\times\UU$
defined by \begin{align}\label{transfo:eq}
(X,U)=\bigl(\varphi_g(x),\psi_g(u)\bigr),
\end{align}
where $\varphi_g$ and $\psi_g$ are local diffeomorphisms. Notice
$\varphi_g$ acts on $\XX$ and $\psi_g$ acts on $\UU$. $u$ can also
denote the time $t$ but in this case $\psi_g$ is the identity
function. The two following definitions are inspired
from~\cite{martin-et-al-cocv03}.

\begin{defn} \label{dyn:inv:def}
The system $\dotex x=f(x,u)$ is \emph{$G$-invariant} if $
f\bigl(\varphi_g(x),\psi_g(u)\bigr)= D\varphi_g(x)\cdot f(x,u)$ for
all $g,x,u$.
\end{defn}
The property also reads $\dotex X = f(X,U)$, i.e., the system
remains unchanged under the transformation \eqref{transfo:eq}.
\begin{defn} \label{output:inv:def}
The output $y=h(x,u)$ is \emph{$G$-equivariant} if there exists a
transformation group $(\varrho_g)_{g\in G}$  on ${\mathcal Y}$ such
that $
h\bigl(\varphi_g(x),\psi_g(u)\bigr)=\varrho_g\bigl(h(x,u)\bigr)$ for
all $g,x,u$.
\end{defn}

With $(X,U)=\bigl(\varphi_g(x),\psi_g(u)\bigr)$ and
$Y=\varrho_g(y)$, the definition means $Y=h(X,U)$. The two previous
definitions can be illustrated by the commutative diagram
\[ \begin{CD}
T{\mathcal X} @>D\varphi_g>>T{\mathcal X}\\
@AfAA @AfAA\\
{\mathcal X}\times{\mathcal U} @>\varphi_g\times\psi_g>> {\mathcal X}\times{\mathcal U}\\
@VhVV @VhVV\\
{\mathcal Y} @>\varrho_g>> {\mathcal Y}
\end{CD} \]

\subsection{Basic assumptions}
From now on we consider a $G$-invariant system $\dotex x=f(x,u)$
 with a $G$-equivariant output
$y=h(x,u)$. We let $r\leq n $ be the dimension of the group~$G$. We
systematically assume for each $x$, the mapping $g \mapsto
\varphi_g(x)$ is full rank.

\subsection{The moving frame method, invariant vector fields, base and fiber
coordinates}\label{preliminaries:sec}

\subsubsection{Moving frame method}\label{mov:frame:sec}
This paragraph is independent of the rest of the paper. It is a
recap of the general presentation
of~\cite[theorem~$8.25$]{olver-book99}). Take a r-dimensional
transformation group $G$ acting
 on~$\Sigma\subset\Rset^{s}$ via the diffeomorphisms $(\phi_g)_{g\in G}$  such that $r\leq s$.
We suppose that $\partial_g\phi_g$ has full rank $r:=\dim G$ at the
point $(e,\xi^0)\in G\times\Sigma$. We can then split $\phi_g$ into
$(\phi_g^a,\phi_g^b)$ with respectively $r$ and $s-r$ components so
that $\phi_g^a$ is invertible with respect to~$g$
around~$(e,\xi^0)$. The \emph{normalization equations} are obtained
 setting $$ \phi^a_g(\xi)=c,$$ with $c$ a constant in the range
of~$\phi^a$. The implicit function theorem ensures the existence of
the local solution $g=\gamma(\xi)$ (the map
$\gamma:\Sigma\rightarrow G$ is known as the \emph{moving frame}).
Thus $$\phi^a_{\gamma(\xi)}(\xi)=c$$One can also say
$\{\phi_e^a(\xi)=c\}$ defines a coordinate cross-section to the
orbits, and $g=\gamma(\xi)$ is the unique group element that maps
$\xi$ to the cross-section. Finally, we get a complete set~$J$ of
$s-r$ functionally independent invariants by substituting
$g=\gamma(\xi)$ into the remaining transformation rules,
\[ J(\xi):=\phi^b_{\gamma(\xi)}(\xi). \]
The invariance property means $J\bigl(\phi_g(\xi)\bigr)=J(\xi)$ for
all~$g,\xi$. To prove it let $\zeta=\phi_g(\xi)$. We have
$\phi^a_{\gamma(\phi_g(\xi))}(\phi_g(\xi))=\phi^a_{\gamma(\zeta)}(\zeta)=c$.
But the group composition implies
$\phi^a_{\gamma(\phi_g(\xi))}(\phi_g(\xi))=\phi^a_{\gamma(\phi_g(\xi))g}(\xi)$.
Thus $\phi^a_{\gamma(\phi_g(\xi))g}(\xi)=c$ which proves by unicity
of $\gamma(\xi)\in G$\begin{align}\label{moving:frame:property:eq}
\gamma(\phi_g(\xi))g=\gamma(\xi)
\end{align}
which is the main property (equivariance) of the moving frame  that
proves indeed
$$J\bigl(\phi_g(\xi)\bigr)=\phi^b_{\gamma(\phi_g(\xi))}(\phi_g(\xi))=\phi^b_{\gamma(\phi_g(\xi))g}(\xi)=\phi^b_{\gamma(\xi)}(\xi)=J(\xi)$$
Moreover any other local invariant $J'$, i.e, any real-valued
function $J'$ which verifies $J'(\phi_g(\xi))=J'(\xi)$ for all
$g,\xi$ can be written as a function of the complete set of
invariants: $J'=\mathcal{H}(J)$.

\subsubsection{Invariant vector fields and invariant
frame}\label{inv:frame:sec}

The moving frame method allows us to build invariant frames, which
play a role in the construction of invariant observers.
\begin{defn}\label{inv:vec}
A vector field $w$ on $\XX$ is said to be G-invariant if the system
$\frac{d}{dt}x=w(x)$ is invariant. This means
$w(\varphi_g(x))=D\varphi_g(x)\cdot w(x)$ for all $g$, $x$.
\end{defn}
\begin{defn}\label{inv:frame}
An invariant frame $(w_1,...,w_n)$ on $\XX$ is a set of n linearly
point-wise independent G-invariant vector fields, i.e
$(w_1(x),...,w_n(x))$ is a basis of the tangent space to $\XX$ at
$x$.
\end{defn}
We are now going to explain how to build an invariant frame. We
follow \cite{olver-book95}, theorem~$2.84$ and we apply the moving
frame method to the following case: $\Sigma=\XX$, and
$(\phi_g)_{g\in G}=(\varphi_g)_{g\in G} $  and the normalization
equations $\varphi_g^a( x)=c$ give $g=\gamma(x)$.
\begin{lem}\label{inv:frame:lem}The vector fields defined by
 \begin{equation}\label{frame:eq}
   w_i(x):= \left(D\varphi_{\gamma(x)}(x)\right)^{-1}\cdot
   \dv{}{x_i},\quad i=1,\dots,n,
 \end{equation}
 where
$(\dv{}{x_1},...,\dv{}{x_n})$ is the canonical frame of $\XX$, form
an invariant frame\footnote{One could take any basis $(e_1,...,e_n)$
of $\XX$ instead of the canonical frame}.\end{lem}
\begin{proof}
They are clearly point-wise linearly independent. Each $w_i$ is
invariant because for any group element $b$ we have
\begin{itemize}
    \item $w_i(\varphi_b(x))=
    (D\varphi_{\gamma(\varphi_b(x))}(\varphi_b(x)))
       ^{-1}\ \dv{}{x_i}$ and thus
       $$(D\varphi_b(x))^{-1}\ w_i(\varphi_b(x))
=\left[D\varphi_{\gamma(\varphi_b(x))}(\varphi_b(x))\
     D\varphi_b(x)\right]^{-1}\ \dv{}{x_i}
     $$

    \item the group structure implies that, for any group
    elements
    $c$, $d$, we have    $
     \varphi_c(\varphi_d)(x)=\varphi_{cd}(x)$; thus
     $$
      D\varphi_{c}(\varphi_d(x))\ D\varphi_d(x)=
     D\varphi_{cd}(x)
     $$
Thus, with $c=\gamma(\varphi_b(x))$ and $d=b$, we have
     $$
D\varphi_{\gamma(\varphi_b(x))}(\varphi_b(x))
      D\varphi_b(x)
      = D\varphi_{\gamma(\varphi_b(x))b}(x)
      ;
$$
\item since $\gamma(\varphi_b(x))b \equiv\gamma(x)$ (eq \eqref{moving:frame:property:eq}), we have
(corresponding to definition \ref{inv:vec})
$$
(D\varphi_b(x))^{-1}\ w_i(\varphi_b(x))=
    (D\varphi_{\gamma(x)}(x))
      ^{-1}\ \dv{}{x_i} = w_i(x)
      .
$$
\end{itemize}
\end{proof}
\subsubsection{Base and fiber coordinates}\label{base:fiber:sec}
We introduce base and fiber coordinates which are useful local
coordinates to express G-invariant systems of definition
\ref{dyn:inv:def}. We suppose from now on that G is the
$r$-dimensional ($r\leq n$) group acting on $\XX\times\UU$ (see
\eqref{transfo:eq}) and for each $x$, the mapping $g \mapsto
\varphi_g(x)$ is full rank. The moving frame method provides a set
of fundamental local invariants $z_b\in\mathbb{R}^{n-r}$  of the
group action on $\XX$ alone. Complete it with $z_a\in\mathbb{R}^{r}$
so that ($z_a,z_b$) form coordinates of $\XX$. These coordinates are
called fiber ($z_a$) and base ($z_b$) coordinates
(see~\cite{olver-book99}). One can always choose $z_a$ such that for
any $g\in G$ the group transformation reads
$\varphi_g(z_a,z_b)=(\varpi(z_a),z_b)$ with $g\mapsto \varpi_g(x)$
invertible for all $x\in\mathbb{R}^{r}$. Let $z=(z_a,z_b)$. Let
$\gamma$ be the moving frame which maps $z$ to the coordinate
cross-section $\{z_a=c\}$. The invariant dynamics
\eqref{dyn:inv:def} writes locally in the new coordinates:
\begin{multline}\label{Dynbf:eq}
\begin{aligned}
\dotex
z_a&=D\varpi_{\gamma(z)^{-1}}f_a(c,z_b,\psi_{\gamma(z)}(u))\\
\dotex z_b&=f_b(c,z_b,\psi_{\gamma(z)}(u))
\end{aligned}
\end{multline}
since the system is invariant. Example \ref{chem:ex} illustrates the
interest of such coordinates.

\subsection{Characterization of invariant pre-observers}
\begin{defn}[pre-observer]\label{preobs:def}
The  system $
 \dotex{\hat{x}} = F(\hat{x},u,y)
$ is a \emph{pre-observer} of~\eqref{eq:sys}-\eqref{eq:output} if
for all $x,u$ $F\bigl(x,u,h(x,u)\bigr)=f(x,u)$.
\end{defn}
The definition does not deal with convergence; if moreover $\hat
x(t)\rightarrow x(t)$ as $t\rightarrow+\infty$ for every (close)
initial conditions, the pre-observer is an (asymptotic)
\emph{observer}.

\begin{defn}The pre-observer
$\dotex{\hat x}=F(\hat x,u,y)$ is \emph{G-invariant} if for all
$g,\hat x,u,y$,
\[ F\bigl(\varphi_g(\hat x),\psi_g(u),\varrho_g(y)\bigr)
 = D\varphi_g(\hat x)\cdot F(\hat x,u,y). \]
\end{defn}
The property also reads $\dotex{\hat X} = F(\hat X,U,Y)$, with
$X=\varphi_g(x)$, $U=\psi_g(u)$ and $Y=\varrho_g(y)$. This means the
pre-observer remains unchanged under the action of $G$ on each of
the three spaces $\XX$, $\UU$, and $\YY$ via (resp.) $\varphi_g$,
$\psi_g$ and $\varrho_g$. Obviously we call \emph{invariant
observer} an asymptotic G-invariant pre-observer.

The assumption that the output is G-equivariant is motivated by the
following result of ~\cite{bonnabel-rouchon-LN05}: if the
pre-observer $\dotex{\hat x}=F(\hat x,u,y)$  is invariant and if the
rank of $F$ versus $y$ is equal to $\dim(y)$, then, the output map
$y$ is $G$-equivariant in the sense of
definition~\ref{output:inv:def}.

 In
general the ``usual" output error $\hat y-y=h(\hat x,u)-y$ does not
preserve the system geometry, hence it will not yield an invariant
pre-observer. The key idea in order to build an invariant (pre-)
observer is to use, as noticed in~\cite{aghannan-rouchon-cdc02}, an
invariant output error instead of the usual output error.

\begin{defn} \label{inverr:defn}The smooth map $(\hat x,u,y)\mapsto E(\hat
x,u,y)\in\Rset^p$ is an \emph{invariant output error} if
\begin{itemize}
\item the map $y\mapsto E(\hat x,u,y)$ is invertible for all~$\hat
x,u$
\item $E\bigl(\hat x,u,h(\hat x,u)\bigr)=0$ for all $\hat x,u$
\item $E\bigl(\varphi_g(\hat x),\psi_g(u),\varrho_g(y)\bigr)
 =E(\hat x,u,y)$ for all~$\hat x,u,y$
\end{itemize}
\end{defn}
The first and second properties mean $E$ is an ``output error", i.e.
it is zero if and only if~$h(\hat x,u)=y$; the third property, which
also reads $E(\hat X,U,Y)=E(\hat x,u,y)$, expresses invariance.

\begin{thm}\label{invobs:char:thm}
$\dotex{\hat x}=F(\hat x,u,y)$ is a G-invariant pre-observer for the
$G$-invariant system $\dotex x=f(x,u)$
 with $G$-equivariant output
$y=h(x,u)$ if and only if
\[F(\hat x,u,y)=f(\hat x,u)
 +\sum_{i=1}^n\LL_i\bigl(I(\hat x,u),E(\hat x,u,y)\bigr)
 w_i(\hat x), \]
where $E$ is an invariant output error, $(\hat x,u)\mapsto I(\hat
x,u)\in\Rset^{n+m-r}$ is a full-rank invariant function, the
$\LL_i$'s are smooth functions such that for all~$\hat x$,
$\LL_i\bigl(I(\hat x,u),0\bigr)=0$, and $(w_1,...,w_n)$ is an
invariant frame.
\end{thm}

Since each $\LL_i$ is smooth and satisfies $\LL_i(I,0)=0$, we can write $\LL_i(I,E)=\bar\LL_i(I,E)\cdot E$ where $\bar\LL_i(I,E)$ is a $p\times1$ matrix
with entries depending on $(I,E)$. Hence,
\begin{align*}
\sum_{i=1}^n\LL_i(I,E)w_i
 &=\sum_{i=1}^n w_i\bigl(\bar\LL_i(I,E)\cdot E\bigr)\\
 &=\begin{pmatrix}w_1& \cdots& w_n\end{pmatrix}
 \begin{pmatrix}\bar\LL_1(I,E)\\ \vdots\\ \bar\LL_n(I,E)\end{pmatrix}E
\end{align*}
The observer can thus be written as
\begin{align}\label{inv:obs:matrix}
F(\hat x,u,y)=f(\hat x,u)
 +W(\hat x)\bar\LL\bigl(I(\hat x,u),E(\hat x,u,y))E(\hat x,u,y)
\end{align}
 where $W(\hat x)=\bigl(w_1(\hat x),..,w_n(\hat x)\bigr)$ and
 $\bar\LL$ is a $n\times p$ matrix whose entries depend on $(I,E)$.
 The observer can be thought of as a gain-scheduled observer with a $n\times p$ gain matrix $W\cdot\bar\LL$ multiplied by the nonlinear error~$E$.

 Notice the theorem says nothing about convergence but only deals
with the structure of the pre-observer.

To prove  theorem \ref{invobs:char:thm} we first prove the following theorem which ensures the existence of a (local) invariant output error. The proof is
constructive and relies on the Cartan moving frame method (see section \ref{mov:frame:sec}).
\begin{thm}\label{thm:ioe} We have the three following statements
\begin{itemize}
 \item there is an invariant output error~$(\hat x,u,y)\mapsto E(\hat x,u,y)$
 \item there is a full-rank invariant function $(\hat x,u)\mapsto I(\hat x,u)\in\Rset^{n+m-r}$
 (a complete set of n+m-r independent scalar invariants)
 \item every other invariant output error reads
\[\tilde E(\hat x,u,y)=\LL\bigl(I(\hat x,u),E(\hat x,u,y)\bigr). \]
where $\LL$ is any smooth function such that $\LL(I,0)=0$ and
$E\mapsto \LL(I,E)$ is invertible.
\end{itemize}
\end{thm}
\begin{proof}
We apply the moving frame method (section \ref{mov:frame:sec}) to
the following case: $\Sigma=\XX\times\UU\times\YY$, and $\phi_g$ is
the composite transformation
\[ \phi_g(\hat x,u,y)
:=\bigl(\varphi_g(\hat x),\psi_g(u),\varrho_g(y)\bigr). \]

Since the action of $G$ on $\XX$ is full rank we can split $\hat
x\mapsto\varphi_g(\hat x)$ into $\varphi_g^a(\hat x)\in\RR^r$, which
is invertible with respect to~$g$, and the remaining part
$\varphi_g^b(\hat x)\in\RR^{n-r}$. The $r$ normalization equations
 \begin{equation}\label{norm:eq}
\varphi_g^a(\hat x)=c
\end{equation}
 can then be
solved and give $g=\gamma(\hat x)$, which can be substituted into
the remaining equations to yield the complete set of $n+m+p-r$
functionally independent invariants
\begin{align}\label{I:eq1}
 I(\hat x,u) &:=\bigl(\varphi^b_{\gamma(\hat x)}(\hat x),
 \psi_{\gamma(\hat x)}(u)\bigr)\\\label{I:eq2}
 J_h(\hat x,y) &:=\varrho_{\gamma(\hat x)}(y).
\end{align}
An invariant output error is then given by
\begin{align}\label{IOE:eq}
E(\hat x,u,y):=J_h\bigl(\hat x,h(\hat x,u)\bigr)-J_h(\hat x,y)
\end{align}

Actually, since it is an invariant function of $\hat x, u$ and $y$,
every invariant output error~$\tilde E$ must have the form
\begin{align*}
 \tilde E(\hat x,u,y) &={\cal F}\bigl(I(\hat x,u),J_h(\hat x,y)\bigr)\\
 &={\cal F}\bigl(I(\hat x,u),J_h\bigl(\hat x,h(\hat x,u)\bigr)-E(\hat
 x,u,y)\bigr)\\
 &=\LL\bigl(I(\hat x,u),E(\hat x,u,y)\bigr)
\end{align*}
We used the fact that $J_h\bigl(\hat x,h(\hat x,u)\bigr)$, which is
by construction invariant, must be a function of~$I(\hat x,u)$
(fundamental invariants of $\hat x$ and $u$).
\end{proof}
We are now able to give the proof of theorem \ref{invobs:char:thm}:

\begin{proof}
The vector field~$F$ in the theorem clearly is a pre-observer. Indeed,
\begin{align*}
 F\bigl(x,u,h(x)\bigr)&=f(x,u)+\sum_{i=1}^n
 \LL_i\Bigl(I( x,u),E\bigl(x,u,h(x)\bigr)\Bigr)w_i(x)\\
 &=f(x,u)+\sum_{i=1}^n\LL_i\bigl(I( x,u),0\bigr)w_i(x)\\
 &=f(x,u)
\end{align*}
By construction, it is invariant.

Conversely, assume $\dotex{\hat x}=F(\hat x,u,y)$ is a G-invariant
observer. It can be decomposed on the point-wise independent~$w_i$'s
as
\[ F(\hat x,u,y)= \sum_{i=1}^{n} F_i(\hat x,u,y) w_i(\hat x),\]
where the $F_i$'s are smooth functions. Since it is a pre-observer,
\begin{align*}
 f(x,u)=F\bigl(x,u,h(x,u)\bigr)
 =\sum_{i=1}^nF_i\bigl(x,u,h(x,u)\bigr)w_i(x).
\end{align*}
Since it is a G-invariant pre-observer
\begin{align*}
 &\sum_{i=1}^nF_i\bigl(\varphi_g(\hat x),\psi_g(u),\varrho_g(y)\bigr)
 ~w_i(\varphi_g(\hat x))
 =D\varphi_g(\hat x)\cdot\sum_{i=1}^n F_i(\hat x,u,y)w_i(\hat x)
\end{align*}
but the $w_i$'s verify $D\varphi_g(\hat x)\cdot w_i(\hat
x)=w_i(\varphi_g(x))$, hence
\[ F_i\bigl(\varphi_g(\hat x),\psi_g(u),\varrho_g(y)\bigr)
 =F_i(\hat x,u,y), \quad i=1,\ldots n.\]
Therefore,
\begin{align*}
 F(\hat x,u,y) &=f(\hat x,u)
 +\bigl[F(\hat x,u,y)-f(\hat x,u)\bigr]\\
 &=f(\hat x,u)+\sum_{i=1}^{n}\Bigl[F_i(\hat x,u,y)
 -F_i\bigl(\hat x,u,h(\hat x,u)\bigr)\Bigr] w_i(\hat x).
\end{align*}
The functions $F_i(\hat x,u,y)-F_i\bigl(\hat x,u,h(\hat x,u)\bigr)$
are clearly invariant; hence by theorem~\ref{thm:ioe}, $F_i(\hat
x,u,y)-F_i\bigl(\hat x,u,h(\hat x,u)\bigr)
 =\LL_i\bigl(I(\hat x,u),E(\hat x,u,y)\bigr)$.
\end{proof}
\subsection{Invariant pre-observer: a constructive method}\label{method:sec}The system must
 be invariant (i.e, unchanged by transformation
\eqref{transfo:eq}) with equivariant output (definition
\ref{output:inv:def}). Thanks to the last theoretic section we can
build all symmetry-preserving pre-observers: \emph{a)} Solve the
normalization equations \eqref{norm:eq}. Build an invariant error
thanks to \eqref{IOE:eq}, and a complete set of scalar invariants
$I$ thanks to \eqref{I:eq1}. \emph{b)} Build an invariant frame
thanks to \eqref{frame:eq}. \emph{c)} The general form of invariant
pre-observers is given by theorem \ref{invobs:char:thm}. There is a
convenient alternative form \eqref{inv:obs:matrix}.
\subsection{Local convergence around an
equilibrium}\label{equi:ssec}

In this paragraph we will show it is always possible to turn an
asymptotic observer with a local gain design into an invariant one
with the same local behavior (see the chemical reactor of
section~\ref{chem:ex}). Indeed consider an equilibrium $(\bar x,
\bar u,\bar y)$ characterized by $f(\bar x,\bar u)=0$ and $\bar y =
h(\bar x,\bar u)$. Assume that the linearized system around this
equilibrium is observable. This means that the pair $(A,C)$ is
observable where
$$
A=\dv{f}{x}(\bar x,\bar u), ~B=\dv{f}{u}(\bar x,\bar
u),~C=\dv{h}{x}(\bar x,\bar u),~D=\dv{h}{u}(\bar x,\bar u)
$$
Consider the following locally asymptotic observer
\begin{align}\label{asymp:obs:eq}
\dotex \hat x= f(\hat x, u) + L(\hat y - y)
\end{align}
 where we have chosen the
observer constant gain matrix $L$ such that $A+LC$ is a stable
matrix. In general, such an observer is not invariant. One can build
an invariant observer with the same linear-tangent approximation,
i.e., a locally asymptotic observer of the form
\eqref{inv:obs:matrix}
$$
\dotex \hat x= F(\hat x,u,y)=f(\hat x,u)
 +W(\hat x)\bar\LL\bigl(I(\hat x,u),E(\hat x,u,y)\bigr) E(\hat x,u,y)
 $$
with
\begin{equation}\label{LinTang:cond}
     \dv{F}{\hat x}(\bar x,\bar u,\bar y) =  A+LC, \quad
 \dv{F}{u}(\bar x,\bar u,\bar y) = B+LD, \quad
 \dv{F}{y}(\bar x,\bar u,\bar y) = -L
\end{equation}
Let us suggest a possible choice for $\bar\LL$ in order to satisfy
the above conditions on $\dv{F}{\hat x}$, $\dv{F}{x}$ and
$\dv{F}{y}$ at the equilibrium. Since $E(\hat x, u ,\hat y)\equiv
0$, by differentiation versus $\hat x$, $u$ and $y$, we have at the
equilibrium
$$
\dv{E}{\hat x } = - \dv{E}{y} C, \quad \dv{E}{u } = -\dv{E}{y} D
$$
Let $V$ denote the $p\times p$ square invertible matrix
$V=\dv{E}{y}(\bar x,\bar u,\bar y)$. Take for instance the constant
matrix$$ \bar\LL=-W(\bar x)^{-1} L V^{-1}$$The choice proposed for
$\bar\LL$ is such that the above conditions~\eqref{LinTang:cond} are
fulfilled. We made an invariant observer with same local behavior as
\eqref{asymp:obs:eq}.

\subsection{Invariant state-error and convergence issue}\label{invariant:state:error}

We have no general constructive procedure to design the gain
functions $\LL_i$'s of theorem~\ref{invobs:char:thm} in order to
achieve systematic asymptotic convergence of $\hat x$ towards $x$
for any non-linear system possessing symmetries. Nevertheless the
way the state estimation error is defined plays a key role in
convergence analysis. Instead of the linear state-error $\hat x-x$,
we will rather consider the following invariant state-error
$$ \eta(x,\hat
x)= \varphi_{\gamma(x)}(\hat x) - \varphi_{\gamma(x)}(x)
$$
where $\gamma(x)$ is defined as the solution of~\eqref{norm:eq} with
respect to $g$. Notice it is equivalent to choose $\hat x$ to make
the normalization and consider $\eta(x,\hat x)= \varphi_{\gamma(\hat
x)}( x) - \varphi_{\gamma( \hat x)}(\hat x) $. A remarkable result
is that the error equation only depends on the trajectory via $I$
($n+m-r$ scalar invariants):

\begin{thm} \label{inverrdyn:lem}
The dynamics of the invariant state-error $\eta(\hat x,
x)=\varphi_{\gamma(x)}(\hat x)-\varphi_{\gamma(x)}( x)$ depends only
on $\eta$ and scalar invariants depending on $x$ and $u$:
$$
\dotex \eta = \Upsilon(\eta,I(x,u))
$$
for some smooth function $\Upsilon$ and  where $I(x,u)$ is  defined
in theorem~\ref{thm:ioe}
\end{thm}
\begin{proof}
The error $\eta$ is an invariant: for all $g\in G$ we have
$\eta(\varphi_g(x),\varphi_g(\hat x))=\eta(x,\hat x)$. Thus
$\dotex\eta(\varphi_g(x),\varphi_g(\hat x))=\dotex\eta(x,\hat x)$,
i.e,
\begin{equation}\label{dot:eta:eq}
\begin{aligned}
&\partial_1\eta D\varphi_g(x) f(x,u)+\partial_2\eta D\varphi_g(\hat
x) F(\hat
x,u,h(x,u))\\
=&\partial_1\eta f(x,u)+\partial_2\eta F(\hat x,u,h(x,u))
\end{aligned}
\end{equation}
where $\partial_i$ denotes the partial differential relative to the
$i$-th variable. Let $\sigma(\hat x, x,
 u)=\dotex\eta(x,\hat x)=\partial_1\eta
f(x,u)+\partial_2\eta F(\hat x,u,h(x,u)) $. The equality
\eqref{dot:eta:eq} expresses that $\sigma(\varphi_g(\hat
x),\varphi_g(x),\psi_g(u))=\sigma (\hat x, x,u)$. Since $\hat
x=\varphi_{\gamma(x)^{-1}}\big(\eta+\varphi_{\gamma(x)}(x)\big)$,
$\sigma$ is an invariant function of the variables $(\eta,x,u)$.
Since $\eta$ is an invariant, every invariant function of
$(\eta,x,u)$ (in particular $\dotex\eta$) is a function of $\eta$,
and of a fundamental set of scalar invariants of $x$ and $u$:
$I(x,u)$.
\end{proof}
Such invariant coordinates are not unique. Any invariant function of
$x, \hat x$ and $u$ equal to zero when $\hat x=x$ can be used as an
invariant state-error to analyze convergence. Since it must be
 a function of the  complete set of $2n+m-r$
invariants $\bigl(I(x,u),\varphi_{\gamma(x)}(\hat x)\bigr)$, it must
be a function of $I(x,u)$ and of the invariant state-error
$\eta(x,\hat  x)$:$$ \mathcal{F}\bigl(I(x,u),\eta(x,\hat x)\bigr)
$$where $I(x,u)$ is a complete set of scalar invariant for the action
of $G$ on $\XX\times \UU$, and $\mathcal{F}\bigl(I,0)=0$ for all
$I$. All examples illustrate the interest of such special
coordinates to analyze convergence.

\section{Examples}\label{examples}

 \subsection{The non-holonomic car}\label{car:ex} Consider a non-holonomic car
whose dynamics is the following:
\begin{align}\label{car:eq}
\dotex x =u\cos\theta,\quad
 \dotex y
 =u\sin\theta,\quad \dotex\theta =uv, ~~~~~~~~~~~~~~~~h(x,y,\theta)=(x,y)
\end{align}

where $u$ is the velocity and $v$ is a function of the steering
angle. We suppose the  output is the measurement of the position
$h(x,y,\theta)=(x,y)$ (using a GPS for instance).

The system is independent of the origin and of the orientation of
the frame chosen, i.e.,  it is invariant under the action of
$G=SE(2)$, the group of rotations and translations. We make the
identification $G=\mathbb{R}^2\times S^1$ thus any element of $G$
writes $(x_g,y_g,\theta_g)\in~\mathbb{R}^2\times S^1$. For any
$(x_g,y_g,\theta_g)\in ~G$ the map $\varphi_{(x_g,y_g,\theta_g)}$
corresponds to the action on $G$ on the state space
$\mathbb{R}^2\times S^1$:
\begin{align*}
\varphi_{\left(
  x_g,y_g, \theta_g
\right)}(
  x, y, \theta)
&=\begin{pmatrix}
  x_g\\ y_g\\ \theta_g
\end{pmatrix}\cdot
\begin{pmatrix}
  x\\ y\\ \theta
\end{pmatrix}=\begin{pmatrix}
  x\cos\theta_g-y\sin\theta_g+x_g\\
  x\sin\theta_g+y\cos\theta_g+y_g\\
  \theta+\theta_g\\
\end{pmatrix}\\
\text{and} ~\psi_{\left(
  x_g,y_g,\theta_g
\right)} (u,v) &=\begin{pmatrix}
  u\\ v
\end{pmatrix}
\end{align*}
The dynamics is indeed invariant in the sense of definition
\ref{dyn:inv:def}. Take $(x_g,y_g,\theta_g)\in~G$ and
$(x,y,\theta)\in~\mathbb{R}^2\times S^1$ and $(u,
v)\in\UU=\mathbb{R}^2$. Set
$\varphi_{(x_g,y_g,\theta_g)}(x,y,\theta)=(X,Y,\Theta)$ and
$\psi_{(x_g,y_g,\theta_g)}(u,v)=(U,V)$ (transformation
\eqref{transfo:eq}). The dynamics in the new variables reads the
same:
$$
 \dotex X =U\cos\Theta,\quad
 \dotex Y =V\sin\Theta,\quad
 \dotex\Theta =UV
$$
The output function is equivariant in the sense of definition
\ref{output:inv:def}: for any $x_g$, $y_g$, $\theta_g$, $x$ and $y$
we have $\varrho_{(x_g,y_g,\theta_g)}(x,y)=\begin{pmatrix}
  x\cos\theta_g-y\sin\theta_g+x_g\\
  x\sin\theta_g+y\cos\theta_g+y_g
\end{pmatrix}$. We apply method of section \ref{method:sec} to build
an invariant pre-observer.
\paragraph{Invariant output error}
The normalization equations~\eqref{norm:eq} write with $c=0$:
\begin{align*}
  x\cos\theta_0-y\sin\theta_0+x_0 &=0\\
  x\sin\theta_0+y\cos\theta_0+y_0 &=0\\
  \theta+\theta_0 &=0
\end{align*}
hence
\[
\begin{pmatrix}
  x_0\\ y_0\\ \theta_0
\end{pmatrix}
=\begin{pmatrix}
  -x\cos\theta-y\sin\theta\\ x\sin\theta-y\cos\theta\\ -\theta
\end{pmatrix}
=\gamma
\begin{pmatrix}
  x\\ y\\ \theta
\end{pmatrix}
\]
A complete set of invariants is given by (see~\ref{I:eq1}):
$I(x,y,\theta,u,v)= \psi_{\gamma(x,y,\theta)}(u,v)
 =\begin{pmatrix}
  u\\ v
\end{pmatrix}$. Let 
$(x_0, y_0, \theta)^T=\gamma(\hat x,\hat y,\hat \theta)^T$. An
invariant output error writes (see~\eqref{IOE:eq}):
\begin{align*}
    E &=\varrho_{(x_0,y_0,\theta_0)}(\hat x,\hat y)
    -\varrho_{(x_0,y_0,\theta_0)}(x,y)\\
    &=\begin{pmatrix}
    \cos\theta_0& -\sin\theta_0\\ \sin\theta_0& \cos\theta_0
    \end{pmatrix}
    \begin{pmatrix}\hat x\\ \hat y\end{pmatrix}
    +\begin{pmatrix}x_0\\ y_0\end{pmatrix}
    -\begin{pmatrix}
    \cos\theta_0& -\sin\theta_0\\ \sin\theta_0& \cos\theta_0
    \end{pmatrix}
    \begin{pmatrix}x\\ y\end{pmatrix}
    -\begin{pmatrix}x_0\\ y_0\end{pmatrix}\\
    &=\begin{pmatrix}
    \cos\hat\theta& \sin\hat\theta\\ -\sin\hat\theta& \cos\hat\theta
    \end{pmatrix}
    \begin{pmatrix}\hat x-x\\ \hat y-y\end{pmatrix}
\end{align*}
\paragraph{Invariant frame} To build an invariant frame we apply formula \eqref{frame:eq}.
Since
$\bigl(D\varphi_{\gamma(x,y,\theta)}(x,y,\theta)\bigr)^{-1}=D\varphi_{\gamma^{-1}(x,y,\theta)}(x,y,\theta)$
and here $\gamma^{-1}(x,y,\theta)=(x,y,\theta)$ an invariant frame
$(w_1,w_2,w_3)$ is given by the image of the canonical basis of
$\mathbb{R}^2\times \mathbb{S}^1$ by $D\varphi_{(x,y,\theta)}$ ,
i.e, the columns of the matrix
\begin{align*} D\varphi_{(x,y,\theta)}(x,y,\theta)\begin{pmatrix}1 \ 0 \ 0 \\ 0 \ 1\
0\\0\ 0\ 1
\end{pmatrix}=\begin{pmatrix}
    \cos\theta & -\sin\theta & 0 \\
    \sin\theta & \cos\theta & 0 \\
    0 & 0 & 1 \\
    \end{pmatrix}=\begin{pmatrix}
    w_1\ w_2\ w_3
    \end{pmatrix}
 \end{align*}
and one can notice it corresponds to the Frenet frame.
\paragraph{Invariant
pre-observer}Any invariant pre-observer reads
(see~\eqref{inv:obs:matrix})
\begin{align}\label{carobs:eq}
\frac{d}{dt} \begin{pmatrix}\hat x \\ \hat y\\
\hat \theta
\end{pmatrix}=\begin{pmatrix} u\cos \hat\theta \\u\sin \hat\theta
\\  uv
\end{pmatrix}+\begin{pmatrix}
    \cos\hat\theta & -\sin\hat\theta & 0 \\
    \sin\hat\theta & \cos\hat\theta & 0 \\
    0 & 0 & 1 \\
    \end{pmatrix}\bar\LL\begin{pmatrix}
    \cos\hat\theta& \sin\hat\theta\\ -\sin\hat\theta& \cos\hat\theta
    \end{pmatrix}
    \begin{pmatrix}\hat x-x\\ \hat y-y\end{pmatrix}
\end{align}
 where $\bar\LL$ is a smooth $3\times 2$ gain matrix whose
entries depend on the invariant error $E$ but also on the invariants
$I(\hat x,\hat y,\hat \theta,u,v)$.
\paragraph{Error equation}
The variable we choose to make the normalization is $(\hat x,\hat
y,\hat \theta)^T$. The invariant state-error thus reads (see
\eqref{invariant:state:error}):
\begin{equation}\label{carobsErr:eq}
\begin{aligned}
\eta&=\gamma(\hat x,\hat y,\hat \theta)\cdot\begin{pmatrix} \hat x\\ \hat y\\
\hat\theta\end{pmatrix}-\gamma(\hat x,\hat y,\hat \theta)\cdot\begin{pmatrix}  x\\  y\\
\theta\end{pmatrix}\\
&=-\begin{pmatrix}\hat x\\
\hat y\\ \hat\theta\end{pmatrix}^{-1}\cdot\begin{pmatrix}  x\\
 y\\\theta\end{pmatrix}=\begin{pmatrix}
    (\hat x-x)\cos\hat\theta+(\hat y-y)\sin\hat\theta\\
    -(\hat x-x)\sin\hat\theta+(\hat y-y)\cos\hat\theta\\
    (\hat\theta-\theta)\\
    \end{pmatrix}
\end{aligned}
\end{equation}
and let us denote by $\eta=(\eta_x,\eta_y,\eta_\theta)^T$ its
coordinates in $\mathbb{R}^2\times \mathbb{S}^1$. Notice the first
two coordinates of the state error coincide with the invariant
output error: $(\eta_x,\eta_y)=(E_x,E_y)$. Direct computations based
on
\begin{itemize}
\item $\begin{pmatrix}\eta_x\\ \eta_y\end{pmatrix}=\begin{pmatrix}
    \cos\hat\theta& \sin\hat\theta\\ -\sin\hat\theta& \cos\hat\theta
    \end{pmatrix}\begin{pmatrix}\hat x-x\\ \hat y-y\end{pmatrix}$
\item  $\dotex\begin{pmatrix}
    \cos\hat\theta& \sin\hat\theta\\ -\sin\hat\theta& \cos\hat\theta
    \end{pmatrix}=(uv+\bar\LL_{31}\eta_x+\bar\LL_{32}\eta_y)\begin{pmatrix}
    \cos(\hat\theta+\pi/2)& \sin(\hat\theta+\pi/2)\\ -\sin(\hat\theta+\pi/2)& \cos(\hat\theta+\pi/2)
    \end{pmatrix}$
    \item $\dotex \begin{pmatrix}\hat x-x\\ \hat y-y\\ \hat \theta -\theta \end{pmatrix}=\small{\begin{pmatrix} u(\cos \hat\theta-\cos\theta) \\u(\sin
    \hat\theta-\sin \theta)
\\ 0
\end{pmatrix}+\begin{pmatrix}
    \cos\hat\theta & -\sin\hat\theta & 0 \\
    \sin\hat\theta & \cos\hat\theta & 0 \\
    0 & 0 & 1 \\
    \end{pmatrix}\bar\LL
    \begin{pmatrix}\eta_x\\ \eta_y\end{pmatrix}}$
\end{itemize}
yield the following autonomous error equation:
\begin{align*}
\frac{d}{dt}\begin{pmatrix}
  \eta_x\\\eta_y\\ \eta_\theta
\end{pmatrix} &= \begin{pmatrix}u(1-\cos \eta_\theta)
+(uv+\bar\LL_{31}\eta_x+\bar\LL_{32}\eta_y)\eta_y\\u\sin
\eta_\theta-(uv+\bar\LL_{31}\eta_x+\bar\LL_{32}\eta_y)
\eta_x\\0\end{pmatrix} +\bar\LL
    \begin{pmatrix}\eta_x\\ \eta_y\end{pmatrix}
\end{align*}
Indeed the invariant error equation is independent of the trajectory
and only depends on the relative quantities $\eta_x$, $\eta_y$ and
$\eta_\theta$ as predicted by theorem \ref{inverrdyn:lem} since here
the invariants $I$ are $(u,v)$.
\paragraph{Convergence of the error system}
We can here tune the gains so that the error system is almost
globally asymptotically convergent. The error equation writes:
\begin{align*}
    \dotex\eta_x
    &=u(1-\cos\eta_\theta) + \bigl(uv+\bar\LL_{31}\eta_x+\bar\LL_{32}\eta_y\bigr)\eta_y +
    \bar\LL_{11}\eta_x+\bar\LL_{12}\eta_y\\
    \dotex\eta_y
    &=u\sin\eta_\theta - \bigl(uv+\bar\LL_{31}\eta_x+\bar\LL_{32}\eta_y\bigr)\eta_x +
    \bar\LL_{21}\eta_x+\bar\LL_{22}\eta_y\\
    \dotex\eta_\theta &=\bar\LL_{31}\eta_x+\bar\LL_{32}\eta_y
\end{align*}
Take
\begin{align} \label{carobsLL:eq}
 \bar\LL=
 \begin{pmatrix} -\abs{u}a& 0\\ 0& -\abs {u}c\\ 0& -ub\end{pmatrix}
 +\begin{pmatrix} 0& ubE_y-uv\\ uv-ubE_y& 0\\ 0& 0\end{pmatrix}
\end{align}
where, $a$, $b$, $c$ are positive scalar constants, reminding
$E_y=\eta_y$ the error equation writes
\begin{equation}\label{er:car:eq}
\begin{aligned}
    \dotex \eta_x
    &=u(1-\cos\eta_\theta) -\abs{u} a\eta_x\\
   \dotex \eta_y
    &=u\sin\eta_\theta - \abs{u} c\eta_y\\
   \dotex\eta_\theta &=-u b\eta_y
   \end{aligned}
\end{equation}
Let us suppose  $\int_{t_0}^\infty |u(t)|dt=+\infty$ for all
$t_0>0$. Consider the regular change of time scale: $ds=|u|dt$, we
have ($\epsilon_1=\pm 1$ is the sign of $u$)
\begin{align*}
    \frac{d}{ds}\eta_x &=\epsilon_1(1-\cos\eta_\theta) -a\eta_x\\
    \frac{d}{ds}\eta_y &=\epsilon_1\sin\eta_\theta -  c\eta_y\\
    \frac{d}{ds} \eta_\theta &=-\epsilon_1 b\eta_y
\end{align*}
with the following triangular structure:
\begin{align*}
    \frac{d^2}{ds^2}\eta_\theta &= -c \frac{d}{ds} \eta_\theta-b \sin
    \eta_\theta
    \\
    \frac{d}{ds}\eta_x &=\epsilon_1(1-\cos\eta_\theta) -a\eta_x
\end{align*}
The first  equation is the dynamics of the damped  non linear
pendulum with the almost globally stable equilibrium
$\eta_\theta=0$. The second equation is just a first order stable
linear system with $\epsilon_1(1-\cos\eta_\theta)$ as source term.
Thanks to the notion of invariant state errors defined
by~\eqref{carobsErr:eq} we proved
\begin{thm}\label{carobs:thm}
Consider the system \eqref{car:eq}. Assume $\int_{t_0}^\infty
|u(t)|dt=+\infty$ for all $t_0>0$. The non-linear observer
\begin{align*}
\frac{d}{dt} \begin{pmatrix}\hat x \\ \hat y\\
\hat \theta
\end{pmatrix}=\begin{pmatrix} u\cos \hat\theta \\u\sin \hat\theta
\\  uv
\end{pmatrix}+\begin{pmatrix}
    \cos\hat\theta & -\sin\hat\theta & 0 \\
    \sin\hat\theta & \cos\hat\theta & 0 \\
    0 & 0 & 1 \\
    \end{pmatrix}\bar\LL\begin{pmatrix}
    \cos\hat\theta& \sin\hat\theta\\ -\sin\hat\theta& \cos\hat\theta
    \end{pmatrix}
    \begin{pmatrix}\hat x-x\\ \hat y-y\end{pmatrix}
\end{align*}with
\begin{align*}
 \bar\LL=
 \begin{pmatrix} -\abs{u}a& 0\\ 0& -\abs {u}c\\ 0& -ub\end{pmatrix}
 +\begin{pmatrix} 0& ubE_y-uv\\ uv-ubE_y& 0\\ 0& 0\end{pmatrix}
\end{align*}
is almost globally asymptotically convergent.
\end{thm}

\subsection{A chemical reactor}\label{chem:sec}
This example illustrates the various definitions of
section~\ref{th:sec} and the construction of invariant
pre-observers. As an interesting by-product, we show that invariant
pre-observers  always produce positive estimated concentrations. In
theorem~\ref{chemobs:thm}, we propose  a gain design that ensures
global asymptotic stability. The use of base and fiber coordinates
and the notion of invariant error play a crucial role in the
convergence analysis.\label{chem:ex} We consider the classical
exothermic reactor of~\cite{aris-amundson-58}. With slightly
different notations, the dynamics reads
\begin{equation}\label{chem:eq}
\begin{aligned}
 \dotex {X^{in}} &=0 \\
 \dotex{X} &= D(t)(X^{in}-X) - k\exp\left(-\frac{E_A}{RT}\right) X \\
 \dotex T &= D(t)(T^{in}(t)-T) + c\exp\left(-\frac{E_A}{RT}\right) X + v(t) \\
 y &=T
 \end{aligned}
\end{equation}
where $(E_A,R,k,c)$ are positive  and known constant parameters,
$D(t)$, $T^{in}(t)$ and $v(t)$ are known time functions and
$D(t)\geq 0$. The available online measure is  $T$: the temperature
inside the reactor. The parameter $X^{in} >0$, the inlet
composition, is unknown. The reactor composition $X$ is not
measured.

These two differential equations correspond to material and energy
balances. Their structure is independent of the units: the equations
write the same whether  they are written in $mol/l$ or in $kg/l$ for
instance. Let us formalize such independence in terms of invariance.
We just consider a change of material unit corresponding to the
following scaling $X\mapsto g X$ and $X^{in}\mapsto g X^{in}$ with
$g > 0$. The group $G$ is the multiplicative group $\RR_+^\ast$.
Take $x=(X^{in},X,T)$ as state and $u=(c,D(t),T^{in}(t),v(t))$ as
known input. The action on $\XX\times \UU$  is defined for each
$g>0$ via the (linear) transformations
$$
 \left(%
 \begin{array}{c}
 X^{in}\\X\\T
 \end{array}%
 \right)
\mapsto \varphi_g(x) =
\left(%
 \begin{array}{c}
 gX^{in}\\gX\\T
 \end{array}%
 \right)
 ,\quad
  \left(%
 \begin{array}{c}
 c\\D\\T^{in}\\v
 \end{array}%
 \right)
\mapsto \psi_g(u) =
\left(%
 \begin{array}{c}
 c/g\\D\\T^{in}\\v
 \end{array}%
 \right).
$$
The dynamics \eqref{chem:eq} is  invariant in the sense of
definition~\ref{dyn:inv:def}. Since $y=T$ is unchanged by $G$
($\varrho_g(y)\equiv y$ here), it is
 a $G$-equivariant output in the sense of
 definition~\ref{output:inv:def}. We apply method of section \ref{method:sec} to build
an invariant pre-observer.

\paragraph*{a) Invariant output error and complete set of invariants}We choose the
second component of $\varphi_g$ for the normalization and take as
normalizing
  equation~\eqref{norm:eq}: $ gX=1$, i.e. $\gamma(x)= 1/X$.
Then using $\eqref{IOE:eq}$ the  invariant output error is  $E(\hat
x, u,y)=\hat T - y$ and using \eqref{I:eq1} the complete set of
invariant $I$ is made of the remaining components of
$\varphi_{1/\hat X}(\hat x)$ and $\psi_{1/\hat X}(u)$: $I(\hat x, u)
= (\hat X^{in}/\hat X, \hat T,c\hat X, D, T^{in}, v)$.
\paragraph*{b) Invariant frame}
According to~\eqref{frame:eq}, an invariant frame is:
$$
w_1= X^{in}\dv{~}{X^{in}}, \quad w_2=X \dv{~}{X}, \quad
w_3=\dv{~}{T}
$$
where $w_1$ has been multiplied by the scalar invariant $X^{in}/X$.
\paragraph*{c) Invariant pre-observer}
According to theorem~\ref{invobs:char:thm}, invariant pre-observers
have the following structure
\begin{equation}\label{chemobs:eq}
 \left\{
 \begin{aligned}
    \dotex \hat X^{in} & =
   \LL_1\left(I(\hat x,u), \hat T - T \right) \hat X^{in}
    \\
    \dotex{\hat X} & =D(t)(\hat X^{in} - \hat X) - k\exp\left(-\frac{E_A}{R\hat T}\right)
    \hat X +
   \LL_2\left(I(\hat x,u), \hat T -T \right)
     \hat X
    \\
    \dotex{\hat T} & = D(t)(T^{in}(t) - \hat T) + c\exp\left(-\frac{E_A}{R\hat T}\right)
    \hat X +v(t)
       +\LL_3\left(I(\hat x,u), \hat T -T \right)
 \end{aligned}
 \right.
\end{equation}
where the  $\LL_i$'s are smooth scalar functions such that
$\LL_i(I,0)\equiv 0$. Any invariant observer preserves the fact that
$\hat X$ and $\hat X^{in}$ are positive quantities. Indeed the
domain $\{(\hat X^{in},\hat X, T)\in\RR^3~| ~ \hat X^{in}>0, \hat X
>0\}$ is positively invariant for~\eqref{chemobs:eq}, whatever the
choices made for  $\LL_1$, $\LL_2$ and $\LL_3$ ($D(t)\geq 0$).
\paragraph*{d) Convergence around an equilibrium}
 Assume that around a steady-state $(\bar
X^{in},\bar X,\bar T)$ of~\eqref{chem:eq}, we designed the three
constant gains $L_1$, $L_2$, and $L_3$, such that
\begin{align*}
    \dotex \hat X^{in} & =
   L_1 (\hat T - T)
    \\
    \dotex{\hat X} & =D(t)(\hat X^{in} - \hat X) - k\exp\left(-\frac{E_A}{R\hat T}\right)
    \hat X +
  L_2 (\hat T - T)
    \\
    \dotex{\hat T} & = D(t)(T^{in}(t) - \hat T) + c\exp\left(-\frac{E_A}{R\hat T}\right)
    \hat X +v(t)
       +L_3 (\hat T - T)
\end{align*}
is locally convergent around $(\bar X^{in},\bar X,\bar T)$. Then
following the procedure of subsection~\ref{equi:ssec},  we get the
invariant observer
\begin{align*}
    \dotex \hat X^{in} & =
   L_1 (\hat T - T) \frac{\hat X^{in}}{\bar X^{in}}
    \\
    \dotex{\hat X} & =D(t)(\hat X^{in} - \hat X) - k\exp\left(-\frac{E_A}{R\hat T}\right)
    \hat X +
  L_2 (\hat T - T) \frac{\hat X}{\bar X}
    \\
    \dotex{\hat T} & = D(t)(T^{in}(t) - \hat T) + c\exp\left(-\frac{E_A}{R\hat T}\right)
    \hat X +v(t)
       +L_3 (\hat T - T)
\end{align*}
that exhibits identical performances around the steady-state.
Moreover it provides  automatically positive estimations for $X$ and
$X^{in}$, and the performances are independent of the choice of
units.
\paragraph*{e) Invariant error and global convergence of the observer}
As the dimension of $G$ is strictly smaller than the dimension of
$\XX$ it is interesting  to write the dynamics with the base and
fiber coordinates of section \ref{base:fiber:sec} which are globally
defined on  the physical domain $\{(X^{in},X, T)\in\RR^3~| ~
X^{in}>0, X >0\}$. Consider the following change of variable:
$$
\left(
  \begin{array}{c}
  X \\ \ X^{in} \\  T
  \end{array}
\right) \mapsto \left(
  \begin{array}{c}
   \ Z= \log(X)
 \\
 \xi= \log ( X/ X^{in})
 \\

 \ T
  \end{array}
\right)
$$
Indeed $X$ corresponds to fiber coordinate and $X/X^{in}, T $ to
base coordinates. We took the $\log$ of these quantities so that the
computation of time derivatives is easier. The dynamics
\eqref{chem:eq} now writes:
\begin{align*}
 \dotex{ Z} & =D(\exp(-\xi)-1)
    -k\exp\left(-\frac{E_A}{RT}\right)
    \\
    \dotex{ \xi} & =D(\exp(-\xi)-1)
    -k\exp\left(-\frac{E_A}{RT}\right)
    \\
 \dotex{ T} & =  D(T^{in}-T)+c\exp
 Z\exp\left(-\frac{E_A}{RT}\right)+v(t)
\end{align*}
and the invariant observer \eqref{chemobs:eq} writes:
\begin{align*}
 \dotex{ \hat Z} & =D(\exp(-\hat \xi)-1)
    -k\exp\left(-\frac{E_A}{R\hat T}\right)
    +  \LL_2\\
    \dotex{ \hat \xi} & =D(\exp(-\hat \xi)-1)
    -k\exp\left(-\frac{E_A}{R\hat T}\right)+ \LL_2- \LL_1\\
 \dotex{ \hat T} & =  D(T^{in}-\hat T)+c\exp\hat
 Z\exp\left(-\frac{E_A}{R\hat T}\right)+v(t)+ \LL_3
\end{align*}
Consider the following gain design ($\beta>0$ and $\kappa>0$ are two
arbitrary parameters)
\begin{align*}
\LL_2&=-\beta c \exp \hat Z \exp\left(-\frac{E_A}{RT} \right) (\hat
T - T) +k \exp\left(-\frac{E_A}{R\hat
T}\right)-k\exp\left(-\frac{E_A}{R T} \right)
 \\
 \LL_2- \LL_1 &= k \exp\left(-\frac{E_A}{R\hat T}\right)-k\exp\left(-\frac{E_A}{R T} \right)
 \\
 \LL_3&=\left(-\kappa c  \exp\left(-\frac{E_A}{RT}\right) (\hat T-T)
 -c \exp\left(-\frac{E_A}{R\hat T}\right)\right.\\
 &\qquad \left.+c\exp\left(-\frac{E_A}{R T} \right)
 \right)\exp \hat Z+D(\hat T-T)
\end{align*}
The choice of such non-linear gains ensure global asymptotic
stability when there exists $M$ and $\alpha >0$ such that  the
measurements verify  for all $t\geq 0$,  $M\geq X^{in}, D(t), T(t)
\geq \alpha$. It implies (see \eqref{chem:eq}) there exists
$\sigma>0$ such that $t\geq 0$, $M\geq X^{in}, X(t), D(t), T(t) \geq
\sigma$. The design, although specific to the example relies on the
notion of invariant state error (see
subsection~\ref{invariant:state:error}). Since the normalizing
  equation~\eqref{norm:eq} is: $ gX=1$, i.e. $\gamma(  x)= 1/ X$ the
   invariant
   state-error writes in the new variables $\eta=(\tilde Z,\tilde \xi,\tilde T)$ where
\begin{align*}
\tilde Z &= \hat Z-Z =\log(\hat X / X)
 \\
 \tilde \xi &= \hat \xi-\xi =\log (\hat X/\hat X^{in})-\log ( X/ X^{in})
 \\
 \tilde T &= \hat T - T
  \end{align*}
The dynamics of the invariant state error is the following:
\begin{align*}
\dotex{\tilde Z} & =D(\exp(-\hat \xi)-\exp(-\xi))
    - \beta c \exp\left(-\frac{E_A}{RT}+Z\right)\exp\tilde Z ~\tilde T
    \\
    \dotex{\tilde \xi} & =D (\exp(-\hat \xi)-\exp(-\xi))
    \\
    \dotex{\tilde T} & =  c\exp\left(-\frac{E_A}{RT}+Z\right)(\exp \tilde Z -1)
        -\kappa c \exp\left(-\frac{E_A}{RT}+Z\right)\exp\tilde Z~ \tilde T
\end{align*}
Since  $M\geq  D(t)\geq \sigma$,   we have $\lim_{t\mapsto+\infty}
(\hat \xi(t) - \xi(t)) = 0$, which means the dynamics of the system
on the base coordinate $\xi$ converges independently from its
initial value. And the system writes:
\begin{align*}
   \dotex{\tilde Z} & =  -  \beta c  \exp\left(-\frac{E_A}{RT}+Z\right)\exp\tilde Z~\tilde T+\epsilon_1(t)
    \\
    \dotex{\tilde T} & =  c\exp\left(-\frac{E_A}{RT}+Z\right)(\exp\tilde Z -1)
       - \kappa c \exp\left(-\frac{E_A}{RT(t)}+Z\right)\exp\tilde Z~ \tilde T
       .
\end{align*}
where $\epsilon_1(t)=D(\exp(-\hat \xi(t))-\exp(-\xi(t)))$ and we
know that $\lim_{t\mapsto+\infty} \epsilon_1(t)=0$ and
$\int|\epsilon_1(t)|dt<\infty$. Consider the regular  change of time
scale $\tau = \int_0^t c\exp\left(-\frac{E_A}{RT(s)}+Z(s)\right)~
ds$. Then:
\begin{align*}
   \frac{d \tilde Z}{d\tau} & =  - \beta \exp\tilde Z ~\tilde T+\epsilon(t)
    \\
    \frac{d\tilde T}{d\tau} & =  (\exp\tilde Z -1)
       - \kappa \exp \tilde  Z ~ \tilde T
\end{align*}where $\epsilon(t)=cX\exp\left(\frac{E_A}{RT}-Z\right)\epsilon_1(t)$.
Take $ V=\tilde Z +\exp(-\tilde Z) + \frac{\beta}{2} \tilde T^2$ as
Lyapounov function.  $\dotex V\leq |\epsilon(t)|(1+2V)$. Thus $V$ is
bounded and so are the trajectories. Let $(\bar Z(t),\bar T(t))$ be
a trajectory. Take $U(t)=V(\bar Z,\bar
T,t)-\int_t^\infty\epsilon(\tau)(1-\exp(-\bar Z (\tau))d\tau$.
$\dotex U= -\beta\kappa {\bar T}^2\exp\bar Z<0$. A standard
application of Barbalat's lemma shows that $(0,0)$  is globally
asymptotically stable.

Guided by invariance considerations,  we have obtained the
\begin{thm}\label{chemobs:thm}
Consider the system \eqref{chem:eq}. Assume there exist $M$ and
$\alpha
>0$ such that for all $t\geq 0$, $M\geq X^{in}, D(t),  T(t)
\geq \alpha$. Then for any $\beta,\kappa>0$
 the following  non-linear observer:
\begin{align*}
    \dotex{\hat X^{in}} & =
   - \beta  \exp\left(-\frac{E_A}{RT(t)}\right)~ (\hat T - T(t))~
     ~c\hat X ~ \hat X^{in}
    \\
    \dotex{\hat X} & =D(t)(\hat X^{in} - \hat X) - \exp\left(-\frac{E_A}{RT(t)}\right)
     \left( k + \beta (\hat T - T(t)) c\hat X \right)  ~\hat X
    \\
    \dotex{\hat T} & =  \exp\left(-\frac{E_A}{RT(t)}\right)
     \left( 1 - \kappa (\hat T -T(t)) \right)~c\hat X  + D(t)(T^{in}(t) - T(t))  +v(t)
\end{align*}
is globally converging.
\end{thm}

 \subsection{Velocity-aided inertial navigation}\label{nav:sec}
In low-cost navigation systems, the relatively inaccurate gyroscopes and accelerometers are ``aided" by velocity measurements (given by an air-data system
or a Doppler radar) and magnetic sensors. The various measurements are then ``merged" according to the (flat-Earth) motion equations of the aircraft,
usually by a gain-scheduled observer or an extended Kalman filter. The convergence analysis, hence the tuning, of such an observer is far from easy. Using
our theory, we derive in this section a simple invariant observer, which yields an error equation independent of the trajectory of the aircraft. The tuning
of the gains to achieve local convergence around any trajectories is thus straightforward.

Simulations illustrate the good behavior of the observer even in the
presence of noise and sensor biases. They moreover indicate that the
domain of convergence of the observer with respect to the initial
condition should be very large (though we have not investigated the
global behavior).

The derivation of the observer and its implementation are strongly
simplified when the body orientation  is described by a quaternion
of length 1 (rather than by Euler angles or a rotation matrix).

\subsubsection{Quaternions} As in \cite{creamer-jas96}, we use the quaternion
parameterization of SO(3) to derive filters for state estimation.
The quaternions are a non commutative group. Any quaternion $q$ can
be written $q=q^0+q^1 e_1+ q^2 e_2+q^3 e_3$ with
$(q^0,q^1,q^2,q^3)\in\RR^4$, the multiplication $\ast$ is defined by
$$
e_1*e_1=-1,~e_1*e_2=-e_2*e_1=e_3 \text{ with  circular permutations}
$$
and the norm of $q$ is  $\sqrt{(q^0)^2+(q^1)^2+(q^2)^2+(q^3)^2}$.
Any vector $\vec{p}\in\mathbb{R} ^3$ can be identified with the
quaternion $p^1 e_1+ p^2 e_2+p^3 e_3$. We will make this
identification systematically. Then one can associate to any
quaternion whose norm is $1$, a rotation matrix $R_q\in SO(3)$
thanks to the following relation: $q^{-1}*\vec{p}*q=R_q\vec{p}$ for
all $\vec{p}$. The subgroup of quaternions whose norm is $1$ is
denoted by $\mathbb{H}_1$. Conversely, to any rotation $R_q$ of
$SO(3)$ are associated two quaternions  $\pm q$ of length $1$. Thus
although the state space in the example is
$SO(3)\times\mathbb{R}^3$,  we will write the elements of $SO(3)$ as
quaternions whose norm is $1$ (denoted by $\mathbb{H}_1$) and the
vectors of $\mathbb{R}^3$ as quaternions whose first coordinate is
equal to $0$. Numerically, quaternions are easier to manipulate and
compute than matrices in $SO(3)$. The wedge product $v\times \omega$
of vectors of $\mathbb{R}^3$ writes for the associated quaternions:
$(v*\omega-\omega*v)/2.$

\subsubsection{Motion equations} The motion of a flying rigid body
(assuming the Earth is flat and defines an inertial frame) is
described by
\begin{equation}\label{nav:dyn:eq}
\begin{aligned}
    \frac{d}{dt} q &=\frac{1}{2}q\ast\omega\\
    \frac{d}{dt} v &=v\times\omega+q^{-1}\ast\vG\ast q+a\\
    y&=(y_v,y_b)=(v,q^{-1}\ast\vB\ast q)
\end{aligned}
\end{equation}
where $(\omega,a)$ are inputs and
 \begin{itemize}
 \item$q$ is the quaternion of norm $1$ representing the orientation of the
body-fixed frame with respect to the earth-fixed frame. Notice the norm of $q$ is left unchanged by the first equation because $\omega$ is a vector
of~$\RR^3$ (i.e. a quaternion whose first coordinate is~$0$).

\item $\omega=\omega^1 e_1 + \omega^2 e_2 + \omega^3 e_3 $ is
the instantaneous angular velocity vector in the body-fixed frame.
\item $v=v^1 e_1 + v^2 e_2 + v^3 e_3 $ is the velocity vector
of the center of mass in the body-fixed frame
\item $\vG=\vG^1 e_1 +\vG^2 e_2 + \vG^3 e_3  $ is the
gravity vector in the earth-fixed frame.
\item $a=a^1 e_1 + a^2 e_2 + a^3 e_3 $ is the
specific acceleration vector, i.e, the aerodynamics forces divided
by the body mass.
\item $\vB=\vB^1 e_1 + \vB^2 e_2 + \vB^3 e_3 $ is the earth's magnetic field
expressed in the earth-fixed frame. \end{itemize}$\vG$ and $\vB$ are
 constant over the flying area. The first equation
describes the kinematics of the body, the second is Newton's force
law. The measurements are $\omega(t)$, $a(t)$, $v(t)$ and
$q^{-1}(t)\ast\vB\ast q(t)$ (measured by gyroscopes, accelerometers,
air data system or Doppler radar and magnetic sensors). Their
coordinates are known in the body-fixed frame. The goal is to
estimate $q$ and $v$ (i.e make a filter and an estimator for $q$
since it is not measured, and a filter for $v$).

\subsubsection{Invariance of the motion equations} From physical consideration,
the symmetries are associated to the group $SE(3)$ which consists of
rotations and translations in dimension 3. We identify (up to  the
multiplication group $\{-1,+1\}$)  $SE(3)$ and the state space
$\XX={\mathbb H}_1\times \RR^3$. For any $(q_g,v_g)\in G$, the map
$\varphi_{(q_g,v_g)}$ corresponds to the following action of $G$ on
$\XX$:
\[
\varphi_{(q_g,v_g)}(q,v)=
\begin{pmatrix}q_g\\ v_g\end{pmatrix}
\cdot\begin{pmatrix}q\\ v\end{pmatrix}
 =\begin{pmatrix}
  q\ast q_g\\
  q_g^{-1}\ast v\ast q_g+v_g
\end{pmatrix}
\]
Here $u= \left(%
 \begin{array}{c}
 a\\ \omega
 \end{array}%
 \right)$. For any
$(q_g,v_g)\in G$, the map $\psi_{(q_g,v_g)}$ is given by
\begin{align*}
 \psi_{(q_g,v_g)}
 (
  a,\omega
 )
 &=
 \begin{pmatrix}
 q_g^{-1}\ast a\ast q_g-v_g\times(q_g^{-1}\ast \vw\ast q_g)\\
    q_g^{-1}\ast\vw\ast q_g
 \end{pmatrix}
\end{align*}
Let us verify that the dynamics is invariant in the sense of
definition~\ref{dyn:inv:def}. Take $(q_g,v_g)\in G$ and $(q,v)\in G$
and $(a,\omega)\in\UU$.  Set (transformation \eqref{transfo:eq})
$$
 \varphi_{(q_g,v_g)}(q,v) = (Q,V), \quad
 \psi_{(q_g,v_g)}(a,\omega) = (A,\Omega)
 .
$$
\begin{align*} \frac{d}{dt} Q
&=\frac{1}{2}q\ast\omega\ast q_g=\frac{1}{2}q\ast q_g\ast
q_g^{-1}\ast\omega\ast q_g=\frac{1}{2}Q\ast\Omega\\
 \frac{d}{dt} V
&=q_g^{-1}(v\times\omega+q^{-1}\ast\vG\ast q+a)\ast q_g\\
&=(V-v_g)\times\Omega+Q^{-1}\ast\vG\ast Q+A+v_g\times(q_g^{-1}\ast
\vw\ast q_g) \\
&=V\times\Omega+Q^{-1}\ast\vG\ast Q+A
\end{align*}thus the dynamics in
the new variables reads the same: it is invariant in the sense of
definition \ref{dyn:inv:def}. The output function
$$
y=(y_v,y_b)=(v,q^{-1}* \vB* q)
$$
is $G$-equivariant in the sense of definition~\ref{output:inv:def}
with
$$
(Y_v,Y_b)=\varrho_{(q_g,v_g)} (y_v,y_b)= (q_g^{-1}\ast y_v\ast
q_g+v_g~,~q_g^{-1}\ast y_b\ast q_g).
$$

\subsubsection{An invariant pre-observer}
We apply method of section \ref{method:sec} to build an invariant
pre-observer.
\paragraph{Invariant output error and complete set of invariants}
The normalization equations~\eqref{norm:eq} write
\begin{align*}
  q\ast q_0 &=1\\
  q_0^{-1}\ast v\ast q_0+v_0&=0
\end{align*}
(where 1 is the unit quaternion: $1+0e_1+0e_2+0e_3$), hence
\[
\begin{pmatrix}
  q_0\\ v_0
\end{pmatrix}
=\begin{pmatrix}q\\ v\end{pmatrix}^{-1}=\begin{pmatrix}
  q^{-1}\\ -q\ast v\ast q^{-1}
\end{pmatrix}
=\gamma
\begin{pmatrix}
  q\\ v
\end{pmatrix}
\]
Using \eqref{I:eq1}, a complete set of invariants is given by
\[
I(q,v,a,\vw)=\psi_{\gamma(q,v)}
\begin{pmatrix}
  \omega\\ a
\end{pmatrix}
=\begin{pmatrix}
  q\ast\omega\ast q^{-1}\\ q\ast(a+v\times\omega)\ast q^{-1}
\end{pmatrix}
\]
Let $(q_0,v_0)$ be $\gamma(\hat q,\hat v)$. Using \eqref{IOE:eq}, an
invariant output error is given by:
\begin{align*}
    E =\varrho_{(q_0,v_0)}(\hat y)
    -\varrho_{(q_0,v_0)}(y)
&=\begin{pmatrix}
    \hat q\ast(\hat v-v)\ast\hat q^{-1}\\
    \vB-\hat q\ast y_b\ast\hat q^{-1}
    \end{pmatrix}
\end{align*}
\paragraph{Invariant frame}In order to make an invariant frame we must
take a basis of the tangent space to the identity element. The
tangent space to the space of quaternions whose norm is $1$ is the
$3$-dimensional set of all quaternions whose first coordinate is
equal to $0$. Let $e_1$, $e_2$, $e_3$ be the canonical basis of that
space, which can be identified with the canonical basis of $\RR^3$.
We apply formula \eqref{frame:eq} with $x=(q,v)$. Since
$\bigl(D\varphi_{\gamma(q,v)}(q,v)\bigr)^{-1}=D\varphi_{\gamma^{-1}(q,v)}(q,v)$
and here $\gamma^{-1}(q,v)=(q,v)$, an invariant frame is given by
the set of $6$ vector fields whose values in $(q,v)$ are the
following
\begin{align*}
D\varphi_{(q,v)}(q,v)\begin{pmatrix}
    e_i\\0
    \end{pmatrix}_{1\leq i\leq 3}=\begin{pmatrix}
    e_i*q\\0
    \end{pmatrix}_{1\leq i\leq 3},\quad D
    \varphi_{(q,v)}(q,v)\begin{pmatrix}
    0\\e_i
    \end{pmatrix}_{1\leq i\leq3}=\begin{pmatrix}
    0\\q^{-1}\ast e_i\ast q
    \end{pmatrix}_{1\leq i\leq3}
\end{align*}

\paragraph{Invariant pre-observer}
According to theorem \eqref{invobs:char:thm} any invariant
pre-observer reads
    \begin{align*}
    \dotex{\hat q} &=\frac{1}{2}\hat q\ast\omega+
    \sum_{i=1}^3\LL_i^q(I,E)~e_i\ast\hat q\\
    \dotex{\hat v} &=\hat v\times\omega+\hat q^{-1}\ast\vG\ast\hat q+a
    +\sum_{i=1}^3\LL_i^v(I,E)~\hat q^{-1}\ast e_i\ast\hat q,
\end{align*}
where the $\LL_i^q, \LL_i^v$ are smooth functions of $I$ and $E$
such that $\LL_i^q(I,0)=0$ and $\LL_i^v(I,0)=0$. To put it into the
alternative form~\eqref{inv:obs:matrix} we decompose $E$ into $(E_v,
E_b)=(\hat q\ast(\hat v- v)\ast\hat q^{-1},
    \vB-\hat q\ast y_b\ast\hat q^{-1})$ and write
\begin{align*}
 \LL_i^q(I,E_v,E_b)&=\bar\LL_{v,i}^q(I,E_v,E_b)\cdot E_v+\bar\LL_{b,i}^q(I,E_v,E_b)\cdot E_b\\
 \LL_i^v(I,E_v,E_b)&=\bar\LL_{v,i}^v(I,E_v,E_b)\cdot E_v+\bar\LL_{b,i}^v(I,E_v,E_b)\cdot E_b,
\end{align*}
where the $\bar\LL_{v,i}^q, \bar\LL_{b,i}^q, \bar\LL_{v,i}^v, \bar\LL_{b,i}^v$ are $3\times1$ matrices with entries depending on $(I,E_v,E_b)$. Hence,
\begin{align*}
\sum_{i=1}^3\LL_i^q(I,E)~e_i\ast\hat q
 &=\left(\sum_{i=1}^3 e_i\bigl(\bar\LL_{v,i}^q(I,E_v,E_b)\cdot E_v
 +\bar\LL_{b,i}^q(I,E_v,E_b)\cdot E_b\bigr)\right)\ast\hat q\\
 &=\left(\sum_{i=1}^3 e_i\bigl(\bar\LL_{v,i}^q\cdot E_v
 +\bar\LL_{b,i}^q\cdot E_b\bigr)\right)\ast\hat q\\
 &=\left(\begin{pmatrix}\bar\LL_{v,1}^q\\ \bar\LL_{v,2}^q\\ \bar\LL_{v,3}^q\end{pmatrix}E_v
 +\begin{pmatrix}\bar\LL_{b,1}^q\\ \bar\LL_{b,2}^q\\ \bar\LL_{b,3}^q\end{pmatrix}E_b\right)\ast\hat q
\end{align*}
Proceeding in the same way with the other correction term $\sum_{i=1}^3\LL_i^v(I,E)~\hat q^{-1}\ast e_i\ast\hat q$, the general invariant pre-observer
finally reads
    \begin{align} \label{observer1:nav}
    \dotex{\hat q} &=\frac{1}{2}\hat q\ast\omega+(\bar\LL^q_v E_v+\bar\LL^q_b E_b)\ast\hat q\\
    \label{observer2:nav}
    \dotex{\hat v} &=\hat v\times\omega+\hat q^{-1}\ast\vG\ast\hat q+a
    +\hat q^{-1}\ast(\bar\LL^v_vE_v+\bar\LL^v_bE_b)\ast\hat q
\end{align}
where $\bar \LL^q_v$, $\bar \LL^q_b$, $\bar \LL^v_v$ and $\bar
\LL^v_b$, are $3\times 3$ gain matrices whose entries depend on the
invariant errors $E_v$ and~$E_b$ and also on the invariants $I(\hat
q,\hat v,a,\vw)$.

As a by-product of the geometric structure of the observer, we automatically have the desirable property that the the norm of $\hat q$ is left unchanged
by~\eqref{observer1:nav}, because $\omega$ and $\bar\LL^q_v E_v+\bar\LL^q_b E_b$ are vectors of~$\RR^3$ (i.e. quaternions with a zero first coordinate).

\paragraph{Error equation}
The invariant state-error as defined in
section~\ref{invariant:state:error} reads
$\eta=\varphi_{\gamma(q,v)}(\hat q,\hat v) -
\varphi_{\gamma(q,v)}(q,v)$. One can write $\eta=(\eta_q,\eta_v)$
where $\eta_q= \hat q\ast q^{-1}-1$ and $\eta_v=q\ast(\hat v-v)\ast
q^{-1}$.But here the state space can be looked at as the group SE(3)
itself so we will consider the equivalent state-error
\begin{align*}\eta_q= \hat q\ast q^{-1},\quad
\eta_v=q\ast(\hat v-v)\ast q^{-1} .
\end{align*}
so that $\eta$ is an error in the sense of group multiplication.
Thus a small error corresponds to $(\eta_q,\eta_v)$ close to the
group unit element $(1,0)$. Its time derivative verifies:
\begin{align*}\frac{d}{dt}\eta_q&=(\frac{1}{2}\hat q\ast\omega+(\bar\LL^q_v E_v+\bar\LL^q_b E_b)\ast\hat
q)\ast q^{-1}-\hat
q*(q^{-1}*\frac{1}{2}q*\omega*q^{-1})\\&=0+(\bar\LL^q_v
E_v+\bar\LL^q_b
E_b)\ast\eta_q\\
\frac{d}{dt}\eta_v&=q*((\hat v-v)\times\omega+\hat
q^{-1}\ast\vG\ast\hat q-q^{-1}\ast\vG\ast q
    \\&+\hat q^{-1}\ast(\bar\LL^v_vE_v +\bar\LL^v_bE_b)\ast\hat q)*q^{-1}
    +\frac{1}{2}q*\omega*(\hat
v-v)*q^{-1}-q*(\hat v-v)*\omega*q^{-1}\\
&=q*(\hat v-v)\times\omega*q^{-1}+q*\omega\times(\hat
v-v)*q^{-1}+\eta_q^{-1}*\vG*\eta_q\\&-\vG
+\eta_q^{-1}\ast(\bar\LL^v_vE_v+\bar\LL^v_bE_b)\ast\eta_q\\
&=\eta_q^{-1}*\vG*\eta_q-\vG
+\eta_q^{-1}\ast(\bar\LL^v_vE_v+\bar\LL^v_bE_b)\ast\eta_q \\
\end{align*}
where $E_v=\eta_q\ast\eta_v\ast\eta_q^{-1}$ and $E_b=\vB-\eta_q\ast
\vB\ast\eta_q^{-1}$. Thus the error system is autonomous:
\begin{align} \label{naverr:eq}
\frac{d}{dt}\eta_q&=\left[\bar\LL^q_v
(\eta_q\ast\eta_v\ast\eta_q^{-1})+\bar\LL^q_b (\vB-\eta_q\ast
\vB\ast\eta_q^{-1})\right]\ast\eta_q
  \\ \notag
\frac{d}{dt}\eta_v&=
\eta_q^{-1}\ast\left[\vG+\bar\LL^v_v(\eta_q\ast\eta_v\ast\eta_q^{-1})+\bar\LL^v_b(\vB-\eta_q\ast
\vB\ast\eta_q^{-1})\right]\ast\eta_q -\vG
\end{align}
It does depend neither on the trajectory, nor on the inputs $\vw(t)$
and $a(t)$. In the general case (see theorem \ref{inverrdyn:lem})
$\dotex\eta$ is a function of $\eta$ and $I$. But here it does not
even depend on $I$.

\paragraph{Convergence of the linearized error system}
Let us suppose $\hat q$ and $\hat v$ are close to respectively $q$
and $v$. First order approximations write $\delta E_v=\delta\eta_v$
and $\delta E_b=-\delta\eta_q
*\vB+
\vB*\delta\eta_q=2\vB\times\delta\eta_q$. Thus the linearized error
equation writes:
\begin{align*}\frac{d}{dt}\delta\eta_q&=\bar\LL^q_v \delta\eta_v+2\bar\LL^q_b
(\vB\times\delta\eta_q)\\
\frac{d}{dt}\delta\eta_v&=2\vG\times
\delta\eta_q+\bar\LL^v_v\delta\eta_v+2\bar\LL^v_b(\vB\times\delta\eta_q)\end{align*}
Let us choose\begin{align*}
    \bar\LL^q_v&=\begin{pmatrix}0 & -M_{12}&0\\ M_{21}& 0& 0\\ 0& 0& 0\end{pmatrix}
    &\quad
    \bar\LL^v_v&=-\begin{pmatrix}N_{11}& 0& 0\\ 0& N_{22}& 0\\ 0& 0& N_{33}\end{pmatrix}\\
    \bar\LL^q_b&=\begin{pmatrix}0 & 0& 0\\ 0& 0& 0\\ -\lambda \vB^2& \lambda \vB^1& 0\end{pmatrix}
    &\quad
    \bar\LL^v_b&=\begin{pmatrix}0& 0& 0\\ 0& 0& 0\\ 0& 0& 0\end{pmatrix}\\
\end{align*}
In (Earth-fixed) coordinates,
\[  \delta\eta_q:=\begin{pmatrix}
 0\\
 \delta\eta_q^1\\
 \delta\eta_q^2\\
 \delta\eta_q^3
 \end{pmatrix},
 \quad
 \delta\eta_v:=\begin{pmatrix}
 \delta\eta_v^1\\
 \delta\eta_v^2\\
 \delta\eta_v^3
 \end{pmatrix}
  \quad\text{and}\quad
\vG=\begin{pmatrix}
 0\\ 0\\\vg
 \end{pmatrix},
\]
The matrices were chosen so that the error system decomposes in four
decoupled subsystems:
\begin{itemize}
\item the longitudinal subsystem
\begin{align}\label{nav1:eq}
    \begin{pmatrix}\delta\dotex \eta_q^2\\ \delta\dotex \eta_v^1\end{pmatrix}
    &=\begin{pmatrix}0& M_{21}\\ -2\vg& -N_{11}\end{pmatrix}
    \begin{pmatrix}\delta\eta_q^2\\ \delta\eta_v^1\end{pmatrix}
\end{align}
\item the lateral subsystem
\begin{align} \label{nav2:eq}
    \begin{pmatrix}\delta\dotex \eta_q^1\\ \delta\dotex \eta_v^2\end{pmatrix}
    &=\begin{pmatrix}0& -M_{12}\\ 2\vg& -N_{22}\end{pmatrix}
    \begin{pmatrix}\delta\eta_q^1\\ \delta\eta_v^2\end{pmatrix}
\end{align}
 \item the vertical subsystem
\begin{align} \label{nav3:eq}
 \delta\dotex \eta_v^3 &=-N_{33}\delta\eta_v^3
\end{align}
 \item the heading subsystem
\begin{align} \label{nav4:eq}
 \delta\dotex \eta_q^3 &=\lambda \vB^3(\vB^1\delta\eta_q^1-\vB^2\delta
 \eta_q^2)-\lambda((\vB^1)^2+(\vB^2)^2)\delta\eta_q^3
\end{align}

\end{itemize}
We can freely assign the eigenvalues of each of the subsystems.
Guided by invariance considerations we obtained the following non
trivial result:
\begin{thm} \label{navobs:thm}
Consider the dynamics \eqref{nav:dyn:eq}. The non-linear observer
\begin{equation}\label{navobs:eq}
\begin{aligned}
    \dotex{\hat q} &=\frac{1}{2}\hat q\ast\omega(t)+(\bar\LL^q_v E_v+\bar\LL^q_b E_b)
    \ast\hat q
    \\
    \dotex{\hat v} &=\hat v\times\omega(t)+\hat q^{-1}\ast\vG\ast\hat
    q+a(t)
    +\hat q^{-1}\ast(\bar\LL^v_vE_v+\bar\LL^v_bE_b)\ast\hat q
\end{aligned}
\end{equation}

with
$$
 E_v= \hat q\ast(\hat v-y_v(t))\ast\hat q^{-1}
 ,\quad
    E_b=\vB-\hat q\ast y_b(t)\ast\hat q^{-1}
$$
and with the constant gain matrices $\bar \LL^q_v$, $\bar \LL^v_v$,
$\bar \LL^q_b$ and $\bar \LL^v_b$ chosen such that the linear
systems ~\eqref{nav1:eq}, \eqref{nav2:eq}, \eqref{nav3:eq} and
\eqref{nav4:eq},
 are asymptotically stable, converges locally and exponentially
around any system trajectory. The invariant estimation state error
obeys an autonomous differential equation~\eqref{naverr:eq}. The
convergence  behavior and Lyapounov exponents  are completely
independent of the system trajectory and of the inputs.
\end{thm}
Simulations below indicate that the convergence is far from being
only local. We suspect much stronger stability. We conjecture that
such non-linear invariant observer is almost globally convergent. It
can not be globally convergent because of the  following "spin"
effect: if $(\eta_q=1,\eta_v=0)$ is a locally asymptotically stable
steady state for the invariant error equation~\eqref{naverr:eq},
$(\eta_q=-1,\eta_v=0)$ is also a locally asymptotically stable
steady-state. From a physical point of view this is not important
since $\hat q$ and $-\hat q$ correspond to the same rotation
$R_{\hat q}$ in $SO(3)$.

\subsubsection{Simulations}\label{simu:ssec}

To obtain realistic values of $\vw$, $v$, $\va$ and $q^{-1}*\vB* q$
all expressed in the body frame, we generated a trajectory of a
simplified VTOL-like aircraft. The flight is modeled the following
way: initially $q$ is the unit quaternion. Let $k$ denote the
downwards vertical axis of the body frame (quaternion $e_3$) and $P$
 the position of the center of mass of the body. We suppose
the motion is such that $k$ is always collinear to $\ddot P-\vG$. We
suppose $q$ corresponds to the rotation which maps $\vG$ to $k$ and
whose rotation axis is collinear to $\vG\times k$.

We suppose initially that $P(0)=\dot P(0)=\ddot P(0)=0$. $P(t)$
 follows  a circular trajectory whose radius is 5 meters, parameterized by
the angle $\theta (t)$. The function $t\mapsto P(t)$ is $C^3$ with
\begin{itemize}
\item For $0\leq t\leq t_1$ we have $0\leq \theta\leq  \pi/2$ and $\ddot
\theta (t)=c(1-\cos (2\pi t/t_1))$ where
$c=\frac{1}{t_1^2}\frac{2\pi^3}{2\pi^2+1}$ is chosen such that
$t_1=2~s$.
\item For $t_1\leq t\leq t_2$ we have $\pi/2\leq \theta\leq 3\pi/2$
and $\ddot \theta(t) = 0$ with $t_2=4.15~s$.
\item For $t_2\leq t\leq t_3$ we have $3\pi/2\leq \theta\leq 2\pi$
and $\ddot \theta(t) = -c(1-\cos (2\pi (t-t_2)/t_1)$ with
$t_3=6.15~s$.
\end{itemize} The drone eventually  stops after having followed a
circle. The maximum horizontal acceleration is approximately $10~
ms^{-2}$. Such inverse kinematic  model provides realistic values
for  $\va(t)$, $\vw(t)$, $v(t)$ and $q^{-1}(t)*\vB *q(t)$
corresponding to this trajectory.  We take $\vB=[\frac{1}{\sqrt
2},0,\frac{1}{\sqrt 2}]^T$.

\begin{figure}[htb]
\centering
\includegraphics*[width=.7\textwidth]{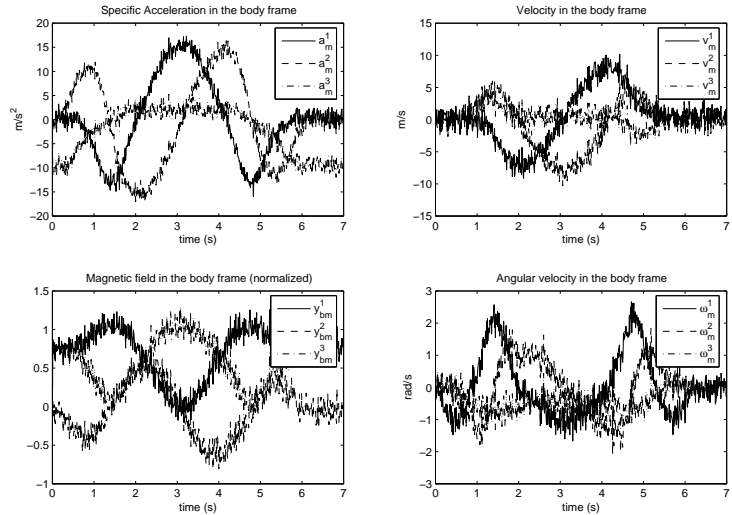}
  \caption{measured signals (with noise  and  bias):   specific acceleration $a$ , velocity  $y_v=v$ ,
  normalized magnetic field $y_b$, and angular velocity $\omega$  in the body frame  }
  \label{m1:fig}
\end{figure}

\begin{figure}[htb]
\centering
\includegraphics*[width=.7\textwidth]{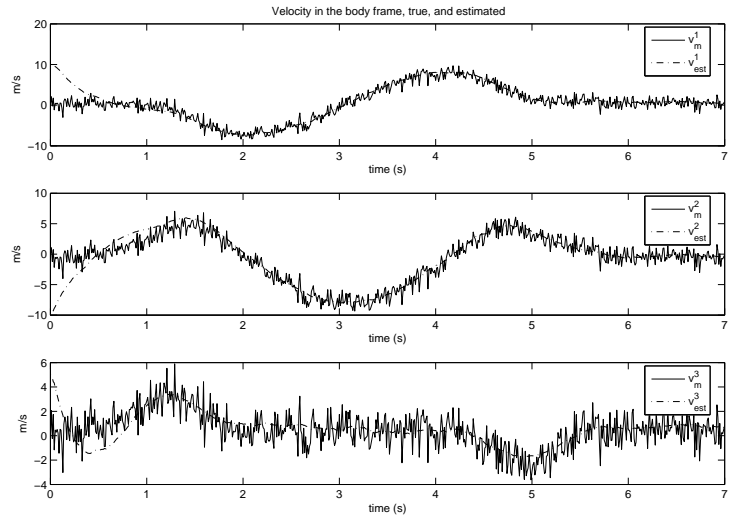}
  \caption{velocity $v$ (solid line)  and estimated velocity  $\hat v$ (dashed-line)
  via the invariant observer~\eqref{navobs:eq} (with noise  and   bias).}
  \label{v1:fig}
\end{figure}

\begin{figure}[htb]
\centering
\includegraphics*[width=.7\textwidth]{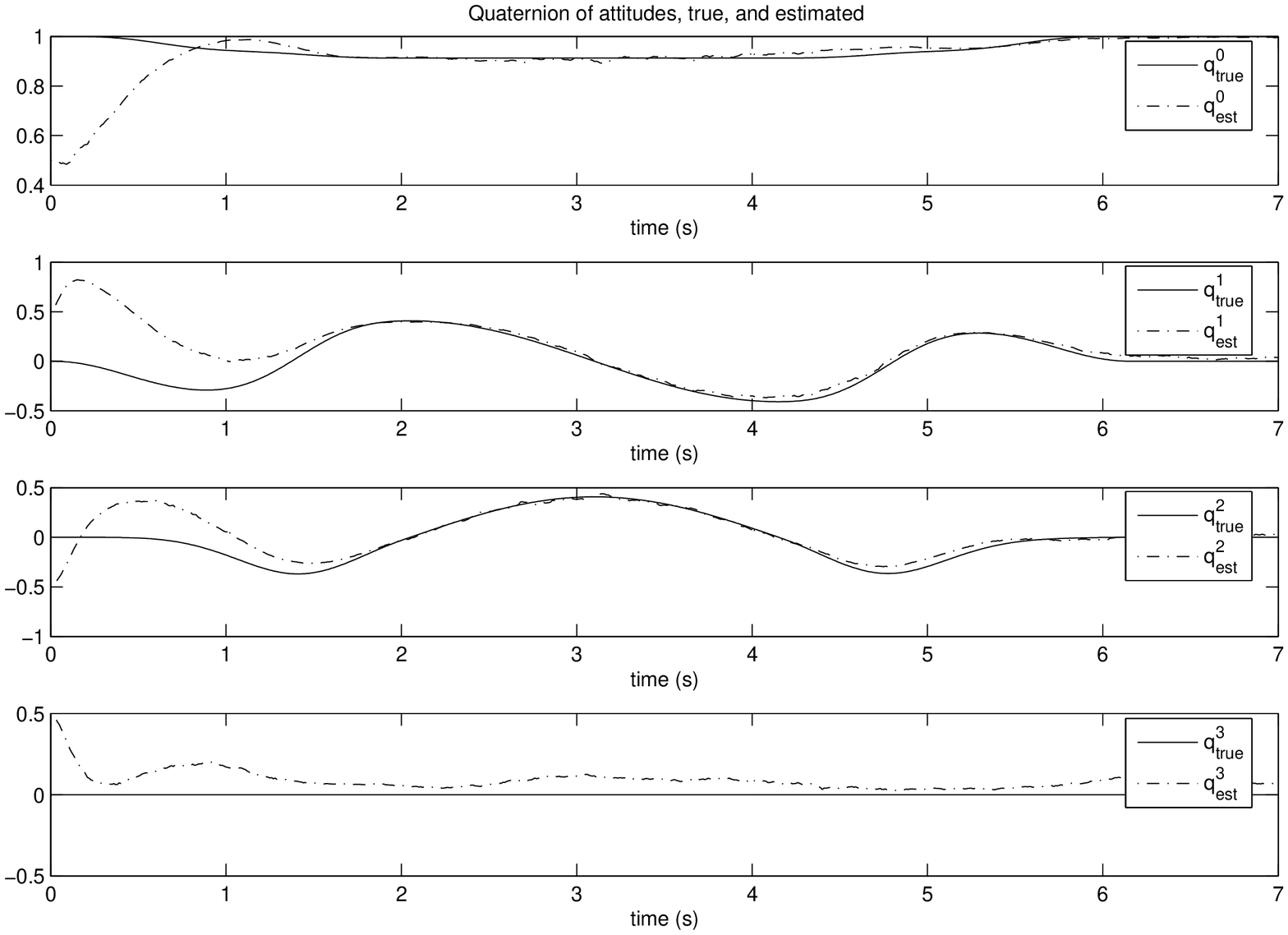}
  \caption{orientation $q$ (solid line) and estimated orientation $\hat q$ (dashed-line)
  via the invariant observer~\eqref{navobs:eq} (with noise  and   bias). }
  \label{q1:fig}
\end{figure}

For the simulations illustrated by
figures~(\ref{m1:fig},\ref{v1:fig},\ref{q1:fig}), the initial
conditions are :
\begin{center}
\begin{tabular}{lcc}
\quad\quad& True system\quad & Observer \eqref{navobs:eq}
 \\
$q_0$ & 1& $\cos (\pi/3)$ \\
$q_1$ & 0& $\sin (\pi/3)/\sqrt(3)$ \\
$q_2$ & 0& $-\sin (\pi/3)/\sqrt(3)$\\
$q_3$ &0 & $\sin (\pi/3)/\sqrt(3)$\\
$v_1$ &0 & 10\\
$v_2$ &0 & -10\\
$v_3$ &0 & 5\\
\end{tabular}
\end{center}
That means the initial rotation differs from the true one up to a
$2\pi/3$ angle. The gains are the following: $M_{12}=M_{21}=0.4$,
$N_{11}=N_{22}=4$, $N_{33}=2$ and $\lambda=4$. With
$\vg=10~ms^{-2}$, the poles of the longitudinal and the lateral
subsystems are: $-2(1\pm i)~s^{-1}$
 and the poles of the vertical and heading subsystems are:
$-2~s^{-1}$.

The  measured signals are noisy and biased: some high frequencies
and some bias are added to the signals $\va(t)$, $\vw(t)$, $v(t)$
and $q^{-1}*\vB*q(t)$ in order to represent the imperfections of the
sensors. The noisy and biased signals are defined by:
$\va_m(t)=\va(t)+0.5~[1,-1,1]+~\sigma_1$, and
$\omega_m(t)=\omega(t)+4\pi/360~[1,-1,1]+0.25~\sigma_2$, and
$v_m(t)=v(t)+0.5~[1,-1,1]+~\sigma_3$, and
$y_{bm}(t)=y_b(t)+0.05~[1,-1,1]+0.1~\sigma_4$, where the $\sigma_i$
are independent normally distributed random 3-dimensional vectors
with mean $0$ and variance $1$, and $\va(t)$, $\vw(t)$, $v(t)$,
$q^{-1}*\vB *q(t)$ are the perfect and smooth signals calculated
from the VTOL-type drone dynamics. These simulations show that the
asymptotic observer~\eqref{navobs:eq} admits a  large attraction
region and is quite robust to measurement noise and bias.

\section{Conclusion}
A theory of symmetry-preserving observers has been developed. It is
mainly composed of: a constructive method to find all the
symmetry-preserving preobservers (see section~\ref{method:sec}), and
a constructive method to find an invariant error between the actual
state of the system and its estimate (see
section~\ref{invariant:state:error}). The resulting invariant error
equation simplifies the convergence analysis. Although we have only
provided examples to support these claims, the following properties
of a symmetry-preserving observer can be expected:
\begin{itemize}
 \item The observer naturally inherits important geometric features of the system (e.g. the observed concentrations in example~\ref{chem:sec} are positive,
 the observed quaternion in example~\ref{nav:sec} has unit norm).
 \item Constant gains can be chosen thanks to the usual linear techniques
(see section \ref{equi:ssec}) to achieve local congergence. If there
are enough symmetries one can expect local convergence around every
trajectory of the system,
 and not only around its equilibrium points or ``slowly-varying"
trajectories.
 \item As the observer respects the geometry of the system, the global behavior tends to be better and the region of attraction larger (compared
 e.g. to a Luenberger observer).
\end{itemize}
Moreover we believe the invariance property of such an observer is
often desirable from an engineering point of view, if not from an
aesthetic one. Finally the method presented in this paper can at
least be seen as a useful new tool in the not-so-full toolbox of
design methods for nonlinear observers, since many physical and
engineering systems exhibit symmetries.

\clearpage


\end{document}